\documentclass[reqno,twoside,11pt,english]{amsart}
\usepackage{amsmath,amsfonts,amssymb,amsthm,epsfig}
\newcommand{\RNum}[1]{\uppercase\expandafter{\romannumeral #1\relax}}
\usepackage{comment}
\usepackage{color}
\usepackage[final,colorlinks,linkcolor=blue,anchorcolor=red,citecolor=blue]{hyperref}
\usepackage{graphicx}
\usepackage{titletoc,epsf}
\usepackage{amsmath,amsfonts,latexsym,amsthm,amsxtra,amssymb,bbm}
\allowdisplaybreaks
\usepackage{graphicx,float}

\newcommand{\lam}{\lambda}
\newcommand{\La}{\Lambda}

\newcommand{\ba}{\backslash}

\newcommand{\la}{\lambda}
\newcommand{\de}{\delta}

\newcommand{\na}{\nabla}
\newcommand{\Om}{\Omega}
\newcommand{\om}{\omega}
\newcommand{\oo}{\infty}

\newcommand{\pa}{\partial}

\newcommand{\va}{\varphi}

\newcommand{\non}{\nonumber}

\newcommand{\al}{\alpha}
\newcommand{\ga}{\gamma}

\newcommand{\dist}{\text{\rm dist}}
\newcommand{\medint}{-\kern -,375cm\int}         
\newcommand{\medintinrigo}{-\kern -,315cm\int}
\newcommand{\wto}{\rightharpoonup}                

\newcommand{\wq}{\infty}
\newcommand{\cal}{\mathcal}

\newcommand{\3}{{\mathcal C}}

\newcommand{\HH}{{\mathcal H}}

\newcommand{\R}{{\mathbb R}}

\newcommand{\N}{{\mathbb N}}

\newcommand{\ep}{\varepsilon}
\newcommand{\su}{\subset}

\newcommand{\ri}{\rightarrow}

\newcommand{\reg}{\text{\rm reg}}
\newcommand{\sing}{\text{\rm sing}}

\newcommand{\Span}{\text{\rm span}}

\numberwithin{equation}{section}
\def\al{\alpha}

\def\ga{\gamma}
\def\de{\delta}
\def\ep{\epsilon}
\def\ka{\kappa}
\def\la{\lambda}

\def\var{\varphi}
\def\om{\omega}
\def\na{\nabla}
\def\Om{\Omega}  
\def\La{\Lambda}

\newtheorem{theorem}{Theorem}[section]

\newtheorem{corollary}[theorem]{Corollary}

\newtheorem{definition}[theorem]{Definition}

\newtheorem{lemma}[theorem]{Lemma}

\newtheorem{proposition}[theorem]{Proposition}
\newtheorem{remark}[theorem]{Remark}
\newtheorem*{remark*}{Remark}

\theoremstyle{definition}

\makeatother

\begin{document}
	
	\title[regularity of fraction harmonic maps]{Optimal first order global regularity for fractional harmonic maps into spheres}

	\author[Y. He, C.-L. Xiang and G.-F. Zheng]{Yu He, Chang-Lin Xiang$^*$ and Gao-Feng Zheng}
	
	\address[Yu He]{School of Mathematics and Statistics, and Key Lab NAA--MOE, Central China Normal University, Wuhan 430079,  P. R.  China}
	\email{yu$\_$he@mails.ccnu.edu.cn}

	\address[Chang-Lin Xiang]{Three Gorges Mathematical Research Center, China Three Gorges University,  443002, Yichang,  P. R. China}
	\email{changlin.xiang@ctgu.edu.cn}

	\address[Gao-Feng Zheng]{School of Mathematics and Statistics, and Key Lab NAA--MOE, Central China Normal University, Wuhan 430079,  P. R.  China}
	\email{gfzheng@ccnu.edu.cn}
	
	\thanks{$^*$Corresponding author: Chang-Lin Xiang}
	
	\date{}

	\begin{abstract}
	This paper is a continuation of our previous work \cite{He-X-Z-24} on global regularity theory for the sphere-valued fractional harmonic maps which was initialized by several authors such as Millot-Pegon-Schikorra \cite{Millot-Pegon-Schikorra-2021-ARMA} and Millot-Sire \cite{Millot-Sire-15}.  By adapting the powerful quantitative stratification approach developed by Naber-Valtorta \cite{Naber-V-2017} on harmonic mappings,  we  establish a corresponding stratification theory for fractional harmonic maps, from which an optimal first order  global regularity theory and the rectifiable structure of singular sets follow.
	\end{abstract}
	
	\maketitle
		
	{\small
		\smallskip
		\newline
		\tableofcontents}
	
	\bigskip
	
		\section{Introduction and main results}

Let  $\Omega\subset\mathbb{R}^n(n\geq 2)$ be a bounded smooth domain and $N\su\mathbb{R}^d$ ($d\ge 2$) be a compact smooth submanifold without boundary. In this paper, we continue the systematic study of the regularity of fractional harmonic maps with spherical targets, building on our previous work \cite{He-X-Z-24} and the foundational contributions of Millot, Pegon, Schikorra and Sire \cite{Millot-Pegon-Schikorra-2021-ARMA, Millot-Sire-15}.

Fractional $s$-harmonic maps are critical points of the nonlocal fractional Dirichlet energy. For a parameter $s\in(0,1)$, the fractional Dirichlet energy functional is defined as
$$\mathcal{E}_s(u,\Om):=
     \frac{\ga_{n,s}}{4}\iint_{(\mathbb{R}^n\times\mathbb{R}^n)
     \ba(\Om^c\times\Om^c)}\frac{|u(x)-u(y)|^2}{|x-y|^{n+2s}}dxdy, $$
where the normalization constant
	$	\gamma_{n, s}=s 2^{2 s} \pi^{-\frac{n}{2}} \Gamma\left(\frac{n+2 s}{2}\right)/{\Gamma(1-s)}	$
	is chosen to ensure consistency with the definition
	$$
	\mathcal{E}_s(u, \Omega)=\frac{1}{2} \int_{\mathbb{R}^n}\left|(-\Delta)^{{s}/{2}} u\right|^2 \mathrm{d} x, \qquad \forall\, u \in \mathcal{D}\left(\Omega ; \mathbb{R}^d\right).
	$$

We consider vector-valued maps $u:\R^n\to N$ belonging to the admissible function space $\widehat{H}^{s}(\Om, N)$, which means that $u$ belongs to the function space
\[\widehat{H}^{s}(\Om, \mathbb{R}^d):=\Big\{u\in L_{loc}^2(\mathbb{R}^n): \mathcal{E}_s(u,\Om)<\oo\Big\}\]
and the pointwise constraint $u(x)\in N$ holds for a.e. $x\in \Om$. To be precise, we adopt the standard definitions of weak, stationary, and minimizing fractional $s$-harmonic maps as follows.
	\begin{definition}\label{uy}(\cite{Millot-Pegon-Schikorra-2021-ARMA})
	We call	$u\in \widehat{H}^{s}(\Om, N)$  a weakly $s$-harmonic map into $N$ if
		\begin{eqnarray}\label{1.1}
			\frac{d}{dt}\Big|_{t=0}\mathcal{E}_s\Big(  \pi_N({u+t\va}),\Om\Big)=0\qquad \text{for all }\,\va\in C_0^\infty(\Omega,\mathbb{R}^d).
		\end{eqnarray}
		Here $\pi_N$ denotes the nearest point orthogonal projection onto $N$.
		
		If a weakly $s$-harmonic map $u\in \widehat{H}^{s}(\Om, N)$ also satisfies
		\begin{eqnarray}\label{1.1}
			\frac{d}{dt}\Big|_{t=0}\mathcal{E}_s\Big(u(x+t\psi(x))\Big)=0\qquad \text{for all }\,\psi\in C_0^\infty(\Omega,\mathbb{R}^n),
		\end{eqnarray}
		then it is called a stationary $s$-harmonic map.
		
		Finally, $u\in \widehat{H}^{s}(\Om, N)$ is  called a minimizing $s$-harmonic map if
		\begin{eqnarray}\label{1.7}
			\mathcal{E}_s(u,\Om)\leq \mathcal{E}_s(v,\Om)
		\end{eqnarray}
		for all $v\in \widehat{H}^{s}(\Omega,N)$ which satisfies $v=u$ on $\mathbb{R}^n\backslash K$ for some compact $K\subset\Omega$.
	\end{definition}

Fractional harmonic maps arise naturally in geometric analysis and mathematical physics, including minimal surfaces with free boundary, nonlocal minimal surfaces, harmonic maps with free boundary value,  Plateau flow problems, and fractional Ginzburg-Landau models for supraconductivity;  see e.g. \cite{Berlyand, A. M. Fraser, Da Lio-book-2018, F. L. T. Riviere, Lio-Pigati-2020-Pisa, Schikorra-2018-Anal, Struwe-1984-Invent, Struwe-2024-APDE} and the references therein. As a canonical nonlocal generalization of classical harmonic maps, the systematic investigation of fractional harmonic maps was initiated by Da Lio and Rivi\`ere  \cite{Da and Ri, Da and Ri 2} for the case $s=1/2$. Subsequent works have extended the theory to general fractional exponents and uncovered deep connections with multiple nonlinear nonlocal problems; we refer to \cite{DaLio-Schikorra-14-ACV, KNS-26-arXiv, Mazowiecka-Schikorra-2018, Millot-Pegon-2020-CV, Millot-Pegon-Schikorra-2021-ARMA, Millot-Sire-15, Millot-Sire-Yu-18, Schikorra-2012-CV, Schikorra-2015-CPDE, Schikorra-2018-Anal} and the references therein for comprehensive developments.

Regularity theory is one of the central issues for fractional harmonic maps.  Early results established smoothness of $1/2$-harmonic maps on the real line (see e.g. \cite{Da and Ri, Da and Ri 2}). In a series of seminal works \cite{Millot-Pegon-2020-CV, Millot-Pegon-Schikorra-2021-ARMA, Millot-Sire-15, Millot-Sire-Yu-18}, Millot, Sire, Pegon, and Schikorra developed a comprehensive partial regularity theory for general $s$-harmonic maps in arbitrary dimensions .  For the most typical and important class of sphere-valued $s$-harmonic maps, their key results can be summarized as follows:

(1) For $n=1$ and $1/2\le s<1$, all weakly $s$-harmonic maps are smooth (\cite[Theorem 1.1]{Millot-Pegon-Schikorra-2021-ARMA});

(2) For stationary $s$-harmonic maps in dimension $n\ge 2$, the singular set
\[\sing(u):=\Om \setminus \reg(u), \qquad \reg(u)=\{x\in \Omega: u \text{ is  continuous in a neighborhood of }x\},\] satisfies the Hausdorff dimension estimate $\dim_{\cal H} \sing(u)\le n-1-[2s]$ (\cite[Theorem 1.2]{Millot-Pegon-Schikorra-2021-ARMA});

(3) For minimizing $s$-harmonic maps, the Hausdorff dimension estimate for singular set is improved to $\dim_{\cal H} \sing(u)\le n-2$ for all $n\ge 2$, with complete regularity $\sing(u)=\emptyset$ in dimension $n=1$, see \cite[Theorem 1.3]{Millot-Pegon-Schikorra-2021-ARMA}. Further regularity improvement can be also found in Millot and Pegon   \cite[Theorem 1.3]{Millot-Pegon-2020-CV}.

For application of regularity theory of $1/2$-harmonic maps, see e.g. \cite{Millot-Sire-15} on $1/2$-Ginzburg-Landau equations and \cite{Lio-Pigati-2020-Pisa} on minimal surfaces with free boundary.
	
While the aforementioned literature provides rigorous partial regularity results (smoothness away from a low-dimensional singular set), the study of global regularity for fractional harmonic maps remained open until our recent contribution \cite{He-X-Z-24}. In that work, we extended the global regularity paradigm of Cheeger-Naber \cite{Cheeger-Naber-2013-CPAM} for classical harmonic maps to the nonlocal setting. For stationary $s$-harmonic maps with $1/2<s<1$ and minimizing $s$-harmonic maps with $0<s\le 1/2$, we established the global gradient integrability $\nabla u\in L^p(D_1)$ for all $1\le p<2$. The key intermediate estimate was an almost-optimal volume bound for tubular neighborhoods of stratified singular sets: for any $\eta>0$, there exists $C=C(n,\Lambda,s,\eta)>0$ such that
    \begin{equation}\label{eq: almost optimal volume estimate}
    \mathrm{Vol}\big(D_r(S_{\eta,r}^k(u)\cap D_1)\big)<Cr^{n-k-\eta}\quad \forall\, 0<r<1,
    \end{equation}
    for all $u\in \widehat{H}_\Lambda^s(D_4,\mathbb{S}^{d-1})$, where the energy-bounded function class is defined by
\begin{equation}\label{eq: bounded energy space}
\widehat{H}^{s}_\La(\Om, N):=\Big\{u\in \widehat{H}^{s}(\Om, N):\mathcal{E}_s(u,\Om)<\La\Big\}
\end{equation}
for fixed constant $\Lambda>0$, and the $r$-tubular neighborhood of $A\subset\R^n$ by
\[ D_r(A)=\{x\in \R^n: \dist(x,A)\le r  \}.\]
As a consequence, the Minkowski dimension of $\sing(u)$ satisfies
$\dim_{\text{Min}}\sing(u)\leq n-2$, which is a slight improvement of the dimension estimate of \cite{Millot-Pegon-Schikorra-2021-ARMA} mentioned in the above.

Despite this progress, the estimate \eqref{eq: almost optimal volume estimate} contains an arbitrarily small but nonvanishing loss exponent $\eta>0$, which prevents one from achieving sharp global regularity and precise geometric characterization of singular sets. In contrast, the breakthrough quantitative stratification theory developed by Naber-Valtorta \cite{Naber-V-2017} yields optimal volume estimates without such error exponents, leading to rectifiability of singular strata and sharp global regularity. This suggests that removing $\eta$-loss in \eqref{eq: almost optimal volume estimate} is the key obstacle to completing the optimal global regularity theory for fractional $s$-harmonic maps. For more applications of the powerful quantitative stratification theory developed by Cheeger, Naber and their group, see e.g. \cite{CJN-21-AnnMath,Guo-Jiang-Xiang-Zheng,Naber-V-2020,Naber-Val-2025-arXiv} and the references therein. 

We emphasize that the borderline case $s=1/2$ exhibits substantially different analytical properties due to distinct compactness behaviors \cite{Millot-Pegon-Schikorra-2021-ARMA, Millot-Sire-15}, and has been treated separately in our recent companion work \cite{GJWXZ-2025-ARMA}. In the present paper, we exclusively consider the generic case
    \[s\in (0,1)\setminus \{1/2\}.\]

Our primary objectives are twofold: (i) establish optimal, loss-free volume estimates for stratified singular sets of spherical-valued stationary $s$-harmonic maps; (ii) deduce the sharp global regularity theory and rectifiable structure of singular sets for fractional harmonic maps via quantitative stratification and symmetry analysis.
We now state our main theorems. The first result provides the optimal volume bounds and rectifiability of all singular strata.

		\begin{theorem}[Optimal volume estimates and rectifiability of singular set]\label{thm80}
		Let  $u\in \widehat{H}^{s}_\La(D_{10}, \mathbb{S}^{d-1})$ be a stationary s-harmonic map. Then, for each $\ep>0$ there exists $C_\ep=C_\ep(n,d,\Lambda,s,\ep)>0$ such that the following volume estimate holds:
		\begin{align}\label{J}
			\mbox{\rm Vol}(D_r(S^k_{\ep,r}(u))\cap D_1)\leq C_\ep r^{n-k},\qquad \forall\,r\in(0,1].
		\end{align}
		Consequently, for all $r\in(0,1]$,
		\begin{align}\label{Jk}
			\mbox{\rm Vol}(D_r(S^k_{\ep}(u))\cap D_1)\leq C_\ep r^{n-k}.
		\end{align}
		Moreover, for each $k$, $S^k_\ep(u)$  are $k$-rectifiable, and also upper Ahlfors $k$-regular. For $\HH^k$-almost every $x\in S^k(u)$, there exists a unique $k$-plane $V^k\subset \R^n$ such that every tangent map of $u$ at $x$ is $k$-symmetric with respect to $V^k$.
\end{theorem}
The symbol $S^k_{\ep,r}(u)$ denotes the $k$-th quantitative stratification of the singular set $\sing(u)$, see Definition \ref{ddd}.	As a direct application of the volume estimate \eqref{J}, we establish the optimal higher regularity for $s$-harmonic maps. We first recall the standard regularity scale function, which was originally introduced by Cheeger-Naber \cite{Cheeger-Naber-2013-CPAM}, adapted to our setting.
	
	\begin{definition}[Regularity scale]\label{pppv}
		Let $u\colon D_2\ri N$ be a map into manifold. For $x\in D_1$, we denote the regularity scale $r_u(x)$ of $u$ at $x$ by
		$$r_u(x):=\max\left\{0\leq r\leq 1:\sup_{D_r(x)}\left(r|\nabla u|\right)\leq1\right\}.$$
	\end{definition}

By definition, one has the set inclusion
$$\big\{x\in D_1:|\nabla u|>r^{-1}\big\}\subset \{x\in D_1:r_u(x)<r\big\},$$
which links the pointwise gradient blow-up to small values of the regularity scale.

	\begin{theorem}[Optimal global regularity estimate]\label{thm84} Assume that $u\in \widehat{H}_\Lambda^s(D_{10},\mathbb{S}^{d-1})$ ($d\ge 2$) is a stationary $s$-harmonic map for $1/2<s<1$, or a minimizing $s$-harmonic map for $0<s< 1/2$. Then there exists a positive constant  $C=C(n,d,\La,s)$ such that
		\begin{equation}\label{emhm}
			\quad {\rm Vol}\Big(\big\{x\in D_1:|\nabla u|>r^{-1}\big\}\Big)\leq{\rm Vol}(\big\{x\in D_1:r_u(x)<r\big\})\leq Cr^2
		\end{equation}
		for all $0<r<1$. Equivalently, the gradient satisfies the optimal global regularity $\na u\in L^{2,\wq}(D_1)$.
	\end{theorem}
	
We provide several remarks concerning the sharpness and necessity of the assumptions in Theorem \ref{thm84}.
\begin{remark}
(i) Sharpness of the regularity exponent. The Lorentz space regularity $\nabla u\in L^{2,\infty}(D_1)$ is optimal in general. A canonical counterexample is the radial projection map $u_{2D}(x)=x/|x|:D^2\to\mathbb{S}^1$, which is stationary for all $s\neq 1/2$ \cite[Remark 1.6]{Millot-Pegon-Schikorra-2021-ARMA} and possesses no better gradient integrability than $L^{2,\infty}$.

(ii) Necessity of the minimizing condition for $0<s<1/2$. For exponents $0<s<1/2$, the stationary condition alone is insufficient to guarantee the global regularity \eqref{emhm}. Indeed, the binary jump map associated to nonlocal minimal surfaces
    \[v(x)=\big(\chi_{\mathbb{R}^n_+}(x)-\chi_{\overline{\mathbb{R}^n_-}}(x),0,\dots,0\big)\]
    is stationary $s$-harmonic for $0<s<1/2$ \cite[Remark 1.5]{Millot-Pegon-Schikorra-2021-ARMA}, \cite[Remark 1.7]{Millot-Sire-Wang-19}, yet its $n$-th weak partial derivative fails to exist near the hyperplane singular set $\mathbb{R}^{n-1}$. Thus the minimizing assumption is essential in the low-exponent regime.

(iii) Comparison with classical harmonic maps. For $1/2<s<1$, our result holds for arbitrary stationary spherical-valued $s$-harmonic maps without additional geometric constraints on the target manifold, in contrast to the classical harmonic map theory in \cite{Naber-V-2017}. This distinction originates from the strong $\widehat{H}^s$-compactness property of stationary fractional $s$-harmonic maps established in \cite{Millot-Pegon-Schikorra-2021-ARMA} for $s\in (0,1)\setminus \{1/2\}$.
\end{remark}

We now outline the main proof strategy of Theorem \ref{thm80} and \ref{thm84} and the contributions of this paper. Our analysis relies on two core tools for nonlocal problems: the Caffarelli-Silvestre extension principle and the monotonicity formula for extended fractional harmonic maps.

First, via the Caffarelli-Silvestre extension \cite{L. Caffarelli and L.Silvestre}, we lift the nonlocal boundary map $u:\mathbb{R}^n\to\mathbb{S}^{d-1}$ to a locally defined map $u^e$ on the upper half-space $\mathbb{R}^{n+1}_+$. A key property is that $u^e$ is smooth in the interior of $\mathbb{R}^{n+1}_+$ and shares exactly the same singular set as the original boundary map $u$. This is quite different from the so-called intrinsic fractional harmonic maps introduced in \cite{M-2011, R-2018}. As a result, this  allows us to study $u^e$ instead of $u$ directly. Moreover, the extended map satisfies a crucial monotonicity formula (see \eqref{1.8}) established in \cite{Millot-Pegon-Schikorra-2021-ARMA}, which enables the adaptation of the quantitative stratification scheme of Naber-Valtorta \cite{Naber-V-2017} to the fractional setting.

To prove Theorem \ref{thm80}, we develop a quantitative symmetry theory for the extended map $u^e$. Unlike classical harmonic maps, the extended problem involves an extra vertical dimension in $\mathbb{R}^{n+1}_+$. We introduce a new notion of boundary quantitative symmetry, combining $0$-homogeneity in the upper half-space with translational invariance restricted to the boundary $\mathbb{R}^n$ (see Definition \ref{quantitative symmetry}). Based on this symmetry notion, we define the quantitative stratification of the singular set and analyze the geometric properties of each stratum through tangent map, including a one-sided Reifenberg regularity property for each stratified singular set (see Proposition \ref{pro32}).

In order to apply the Naber-Valtorta argument, we first verify the Reifenberg condition, a multiscale condition associated with the Jones' number. By virtue of the energy-drop functional $W_r$, we establish refined estimates for the Jones' number. Combining two Reifenberg-type theorems from \cite{Naber-V-2017} and the improved covering argument proposed in \cite{Naber-V-2018}, we derive optimal loss-free volume bounds. We further introduce minor modifications to the covering procedure to simplify the iteration scheme and enhance the clarity of the proof.

Finally, Theorem \ref{thm84} follows by combining the optimal volume estimates with an $\varepsilon$-regularity principle (see Theorem \ref{thm: sym-to-reg}): sufficiently high quantitative symmetry of a stationary $s$-harmonic map implies improved local regularity, which translates to the global Lorentz space gradient bound. See Section \ref{sec: proof of main results} for details.

Before ending this section, we remark on related work for alternative fractional harmonic map models. Moser \cite{M-2011} and Roberts \cite{R-2018, M-R-2022} introduced the class of intrinsic fractional harmonic maps, which do not depend on the imbedding $N\hookrightarrow \mathbb{R}^d$. They proved a partial regularity theory for intrinsic $s$-harmonic maps, see e.g. \cite{M-2011, M-R-2022, R-2018}. In our recent work \cite{W-X-Z-2026}, quantitative stratification of Cheeger-Naber \cite{Cheeger-Naber-2013-CPAM} has been applied to establish almost-optimal global regularity for minimizing intrinsic $s$-harmonic maps.

	\textbf{Notations.} Throughout the paper, we will use the following notations:
	\begin{itemize}
		\item Points in $\mathbb{R}^{n+1}$ are written as $\mathbf{x} = (x,z)$ with $x\in \mathbb{R}^n$ and $z\in\mathbb{R}$;
		\item The boundary $\pa\mathbb{R}^{n+1}_+=\mathbb{R}^n\times\{0\}$ is often identified with $\mathbb{R}^n,$ and subset $A\su\mathbb{R}^n$ is also identified with $A\times\{0\}\subset \pa\mathbb{R}^{n+1}_+.$
		\item   $B_r(\mathbf{x} )$ denotes the
		open ball in $\mathbb{R}^{n+1}$ of radius $r$ centered at $\mathbf{x} = (x,z)$, and $B^+_r(\mathbf{x})=B_r(\mathbf{x} ) \cap \mathbb{R}^{n+1}_+$ denotes the upper half open ball;
		\item $D_r(x)$ denotes the open ball/disk  in $\mathbb{R}^n$ centered at $x$; we write $B_r=B_r(0)$ and $D_r=D_r(0)$ for origin-centered balls;
		\item For any set $A\subset \R^n$, we denote the $r$-tubular neighborhood of $A$ in $\R^n$ by
		\[ D_r(A)=\{x\in \R^n: \dist(x,A)\le r  \};\]
\item For an arbitrary set $G\su\mathbb{R}^{n+1}$, we write
	$$G^+:=G\cap\mathbb{R}^{n+1}_+\ \ \ \ \mbox{and}\ \ \ \ \pa^+G:=\pa G\cap\mathbb{R}^{n+1}_+.$$
	If $G\su\mathbb{R}^{n+1}_+$  is a bounded open set, we shall say that $G$ is admissible whenever 	$\pa G$ is  Lipschitz regular; and  the (relative) open set $\pa^0 G\su\pa\mathbb{R}^{n+1}_+$ defined by
	$$\pa^0G:=\Big\{\mathbf{x}\su\pa G\cap\pa \mathbb{R}^{n+1}_+:B_r^+(\mathbf{x})\su G\ \ \mbox{for some }\ r>0\Big\},$$
	is nonempty and has Lipschitz boundary, such that $\pa G=\pa^+G\cup\overline{\pa^0G}$.\end{itemize}

	\section{Preliminaries}\label{sec:partial regularity}
	
In this section, we introduce  necessary notions and fundamental regularity theory established by e.g. Millot-Pegon-Schikorra \cite{Millot-Pegon-2020-CV, Millot-Pegon-Schikorra-2021-ARMA} for the singular stratification analysis of $s$-harmonic maps.
	
The basic tool we will use is the extension principle of Caffarelli-Silvestre  \cite{L. Caffarelli and L.Silvestre} (see also the probability literature \cite{MO-69}).  First we introduce some function spaces over an open set $G \subseteq \mathbb{R}^{n+1}$. Following \cite{Millot-Pegon-Schikorra-2021-ARMA, Millot-Sire-15}, we introduce the weighted $L^2$-space
	$$
	L^2\left(G,|z|^a \mathrm{d} \mathbf{x}\right):=\left\{v \in L_{\text {loc }}^1(G):|z|^{\frac{a}{2}} v \in L^2(G)\right\}
	$$
	with $a:=1-2 s$	and norm
$\|v\|_{L^2\left(G,|z|^a \mathrm{d} \mathbf{x}\right)}^2:=\int_G|z|^a|v|^2 \mathrm{d} \mathbf{x}$, and then define the weighted Sobolev space
	$$
	H^1\left(G,|z|^a \mathrm{d} \mathbf{x}\right):=\left\{v \in L^2\left(G,|z|^a \mathrm{d} \mathbf{x}\right): \nabla v \in L^2\left(G,|z|^a\mathrm{d} \mathbf{x}\right)\right\},
	$$
equipped with the standard	norm
	$$
	\|v\|_{H^1\left(G,|z|^a \mathrm{d} \mathbf{x}\right)}:=\left(\|v\|^2_{L^2\left(G,|z|^a \mathrm{d} \mathbf{x}\right)}+\|\nabla v\|^2_{L^2\left(G,|z|^a \mathrm{d} \mathbf{x}\right)}\right)^{1/2} .
	$$
	It follows that both $L^2\left(G,|z|^a \mathrm{d} \mathbf{x}\right)$ and $H^1\left(G,|z|^a \mathrm{d} \mathbf{x}\right)$ are separable Hilbert spaces when equipped with the scalar product induced by their respective Hilbertian norms.
Define the weighted Dirichlet energy $\mathbf{E}_s(\cdot, G)$ by
	\begin{equation}\label{E_s}
		\mathbf{E}_s(v, G):=\frac{\boldsymbol{\delta}_s}{2} \int_G|z|^a|\nabla v|^2 \mathrm{d} \mathbf{x} \quad \qquad \text{for } u\in  H^1\left(G,|z|^a \mathrm{d} \mathbf{x}\right),
	\end{equation}
	where $\boldsymbol{\delta}_s:=2^{2 s-1}{\Gamma(s)}/{\Gamma(1-s)}$ is a normalization constant. It was proven that, on a bounded admissible open set $G \subseteq \mathbb{R}_{+}^{n+1}$, the continuous embedding  $L^2\left(G,|z|^a \mathrm{d} \mathbf{x}\right)\hookrightarrow L^\gamma(G)$ holds for every $1\leq\gamma<\frac{1}{1-s}$ whenever $s\in(0,1 / 2)$ by H\"{o}lder's inequality, and $L^2\left(G,|z|^a \mathrm{d} \mathbf{x}\right) \hookrightarrow L^2(G)$ holds for  $s \in[1/2,1)$  since in this case $a \leq 0$. Both cases  imply the continuous embedding
	$$
	H^1\left(G,|z|^a \mathrm{d} \mathbf{x}\right) \hookrightarrow W^{1, \gamma}(G),\qquad \forall\,1<\gamma<\min \{1/(1-s), 2\}.
	$$
Consequently, the following compact embedding holds:
	\begin{equation}\label{eq: comp-embedding}
		H^1\left(G,|z|^a \mathrm{d} \mathbf{x}\right)\hookrightarrow \hookrightarrow  L^\gamma(G), \qquad \forall\,1<\gamma<\min \{1/(1-s), 2\}.
	\end{equation}	
	
The $s$-harmonic extension of a given a measurable function $u:\R^n\to \R$ to the half-space $\mathbb{R}_{+}^{n+1}$ is defined by
	\begin{equation}\label{u^e}
		u^e(x, z):=\sigma_{n,s} \int_{\mathbb{R}^n} \frac{z^{2s} u(y)}{\left(|x-y|^2+z^2\right)^{\frac{n+2s}{2}}} d y,\qquad (x,z)\in \R^n\times (0,\wq),
	\end{equation}
	where $ \sigma_{n, s}:=\pi^{-\frac{n}{2}} {\Gamma\left(\frac{n+2 s}{2}\right)}/{\Gamma(s)} $ is a normalization constant.
	It is well known that $u^e$ satisfies
	\begin{equation}\label{eq: s-harmonic equation}
		\begin{cases}\operatorname{div}\left(z^a \nabla u^{\mathrm{e}}\right)=0 & \text { in } \mathbb{R}_{+}^{n+1}, \\ u^e=u & \text { on } \partial \mathbb{R}_{+}^{n+1}.\end{cases}
	\end{equation}
The following  lemma shows that problem \eqref{eq: s-harmonic equation}  has very good properties for later study on the symmetry of $u$.
\begin{lemma}\label{lem: unique continuation property}(Unique continuation principle)
Suppose $u,v$ solve the first equation of \eqref{eq: s-harmonic equation} and $u\equiv v$ in an open subset of $\R^{n+1}_+$, then $u\equiv v$ holds in $\R^{n+1}_+$. As a result, $u\equiv v$ on $\R^n=\pa \R^{n+1}_+$.
\end{lemma}
	This is simply an application of the unique continuation principle  of \cite[Theorem 1.2]{Garofalo-Lin-86-Indiana}, since \eqref{eq: s-harmonic equation}  is  locally uniformly elliptic in $\R^{n+1}_+$ with smooth coefficients. This property will be used to study symmetry of fractional harmonic maps later.
	
	We will also need the following variational characterization of harmonic extension (see \cite{L. Caffarelli and L.Silvestre}), that is, 	for every   $u \in H^s\left(\mathbb{R}^n\right)$, there holds
	\begin{equation}\label{H^sE^s}
		\begin{aligned}
			{[u]_{H^s\left(\mathbb{R}^n\right)}^2 } & =2\mathbf{E}_s\left(u^{\mathrm{e}}, \mathbb{R}_{+}^{n+1}\right) \\
			& =\inf \left\{2\mathbf{E}_s\left(v, \mathbb{R}_{+}^{n+1}\right): v \in \dot{H}^1\left(\mathbb{R}_{+}^{n+1},|z|^a \mathrm{d} \mathbf{x}\right), v=u \text { on } \mathbb{R}^n\right\},
		\end{aligned}
	\end{equation}
	where the function space $H^s(\R^n)$ is defined as follows: for any open subset $\Om\subset\R^n$, the function space $H^s(\Om)$ consists of all measurable functions $u\in L^2(\Om)$ which satisfies
	\[
	{[u]_{H^s(\Omega)}^2 } :=\frac{\gamma_{n, s}}{2} \iint_{\Omega \times \Omega} \frac{|u(x)-u(y)|^2}{|x-y|^{n+2 s}} \mathrm{~d} x \mathrm{~d} y<\infty. 	 \]
	If $u \in \widehat{H}^s(\Omega)$ for some open set $\Omega \subseteq \mathbb{R}^n$, the following estimates on $u^{\mathrm{e}}$ somehow extends the first equality in \eqref{H^sE^s} to localized setting.	
	\begin{lemma}(\cite[Lemma 2.9]{Millot-Pegon-Schikorra-2021-ARMA}) Let $\Omega \subseteq \mathbb{R}^n$ be an open set. For every $u \in \widehat{H}^s(\Omega)$, the extension $u^{\mathrm{e}}$ given by \eqref{u^e} belongs to $H^1\left(G,|z|^a \mathrm{d} \mathbf{x}\right) \cap L_{\mathrm{loc}}^2\left(\overline{\mathbb{R}_{+}^{n+1}},|z|^a \mathrm{d} \mathbf{x}\right)$ for every bounded admissible open set $G \subseteq \mathbb{R}_{+}^{n+1}$ satisfying $\overline{\partial^0 G} \subseteq \Omega$. In addition, for every point $\mathbf{x}_0=\left(x_0, 0\right) \in \Omega \times\{0\}$ and $r>0$ such that $D_{3 r}\left(x_0\right) \subseteq \Omega$,
		\begin{equation}
			\left\|u^{\mathrm{e}}\right\|_{L^2\left(B_r^{+}\left(\mathbf{x}_0\right),|z|^a \mathrm{d} \mathbf{x}\right)}^2 \leq C\left(r^2 \mathcal{E}_s\left(u, D_{2 r}\left(x_0\right)\right)+r^{2-2 s}\|u\|_{L^2\left(D_{2 r}\left(x_0\right)\right)}^2\right),
		\end{equation}
		and
		\begin{equation}\label{eq: E_sbounded}
			\mathbf{E}_s\left(u^{\mathrm{e}}, B_r^{+}\left(\mathbf{x}_0\right)\right) \leq C_1 \mathcal{E}_s\left(u, D_{2 r}\left(x_0\right)\right),
		\end{equation}
		for a constant $C_1=C_1(n, s)$.
	\end{lemma}
	
	Consequently, there follows	
	
	\begin{corollary} Let $\Omega \subseteq \mathbb{R}^n$ be an open set, and $G \subseteq \mathbb{R}_{+}^{n+1}$ a bounded admissible open set such that $\overline{\partial^0 G} \subseteq \Omega$. The extension operator $u \mapsto u^{\mathrm{e}}$ defines a continuous linear operator from $\widehat{H}^{s}(\Omega)$ into $H^1\left(G,|z|^a \mathrm{d} \mathbf{x}\right)$.
	\end{corollary}

Another natural consequence of the above extension is the simple observation on singularity of functions. We formally write down this observation for  convenience.

\begin{remark}\label{singular set}
	Since $u^e$ is smooth in $\R^{n+1}_+$ by \eqref{u^e}, we have
	$$\sing(u^e)=\sing(u)\times\{0\}.$$
	This allows us to study the singular set of $u$ by examining the singular set of $u^e$.
\end{remark}
	
Now we can state the monotonicity formula of stationary $s$-harmonic maps which plays a central role in the regularity issue.
	\begin{theorem}[Monotonicity formula,~\cite{Millot-Pegon-Schikorra-2021-ARMA,Scheven-2006-MZ}]\label{thm1}
		Let $\Om\su\mathbb{R}^n$ be a bounded open set. If $u\in \widehat{H}^{s}(\Omega,N)$ is stationary in $\Om$, then for every $\mathbf{x}_0=(x_0,0)\in\Om\times\{0\},$  the ``density function''		
		$$r\in(0,dist(x_0,\Om^c))\mapsto\mathbf{\Theta}_s(u^e,\mathbf{x}_0,r)
:=\frac{1}{r^{n-2s}}\mathbf{E}_s(u^e,B_r^+(\mathbf{\mathbf{x}}_0))$$
		is nondecreasing. Moreover,		
		\begin{eqnarray}\label{1.8}			\mathbf{\Theta}_s(u^e,\mathbf{x}_0,r)-\mathbf{\Theta}_s(u^e,\mathbf{x}_0,\rho)= \boldsymbol{\delta}_s\int_{B_r^+(\mathbf{x}_0)\ba B_\rho^+(\mathbf{x}_0)}z^a\frac{|(\mathbf{x}-\mathbf{x}_0)\cdot\na u^e|^2}{|\mathbf{x}-\mathbf{x}_0|^{n+2-2s}}d\mathbf{x}
		\end{eqnarray}	
		for every $0<\rho<r<dist(x_0,\Om^c)$, where $\mathbf{\delta}_s$ is the constant defined in \eqref{E_s}.
	\end{theorem}
	The first important consequences of the  above formula is  the following partial regularity result proved by \cite[Theorem 4.18 and Remark 4.24 ]{Millot-Sire-15} and  \cite[Theorem	5.1]{Millot-Pegon-Schikorra-2021-ARMA}. We state it in the form of scaling function $r_u$ as that of \cite[Proposition 2.5]{He-X-Z-24}.
	
	\begin{theorem}[Partial regularity] \label{thm: partial regularity of MPS}
		There exist $\varepsilon_{1}=\varepsilon_{1}(n,s)>0$ and $\kappa_{2}=\kappa_{2}(n,s)\in(0,1)$
		such that the following holds. Let $u\in\widehat{H}^{s}\left(D_{2};\mathbb{S}^{d-1}\right)$
		be a weakly $s$-harmonic map on $D_{2}$ such that the function
		$r\in(0,2-|\mathbf{x}|)\mapsto\boldsymbol{\Theta}_s\left(u^{\mathrm{e}},\mathbf{x},
		r\right)$
		is nondecreasing for every $\mathbf{x}\in\partial^{0}B_{2}^{+}$.
		If
		\[
		\boldsymbol{\Theta}_s\left(u^{\mathrm{e}},\mathbf{0},1\right)\leqq\varepsilon_{1},
		\]
		then $u\in C^{0,1}\left(D_{\kappa_{2}}\right)$ and
		\[
			r_{u}(0)\ge\kappa_{2}.
		\]
\end{theorem}
	
Combining Marstrand's theorem in geometric measure theory (see e.g. \cite{Mattila-book}), the second important consequence of the  above monotonicity formula is a compactness theory obtained by \cite[Section 7]{Millot-Pegon-Schikorra-2021-ARMA}.
	
	\begin{theorem}\label{thm: compactness} (1) Assume that $s\in (0,1)\backslash\{1/2\}$ and $n>2s$. Let $\{u_{i}\}_{i\ge1}\subset\widehat{H}_{\Lambda}^{s}(D_{4},\mathbb{S}^{d-1})$
		be a sequence of uniformly bounded stationary $s$-harmonic map and
		$u_{i}\wto u$ in $\widehat{H}^{s}(D_{4},\R^d)$ and almost everywhere in $\R^n$. Then $u$ is a stationary
		$s$-harmonic map in $D_{4}$, and for any open subset $\om\subset D_{4}$
		and every bounded admissible open set $G\subset\R_{+}^{n+1}$ satisfying
		$\bar{\om}\subset D_{4}$ and $\overline{\pa^{0}G}\subset D_{4}$,
		there holds
		\[
		\begin{aligned} & u_{i}\to u\qquad\,\,\,\text{strongly in }\widehat{H}^{s}(\om,\R^{d}).\\
			& u_{i}^{e}\to u^{e}\qquad\text{strongly in }H^{1}(G;\R^{d},|z|^{a}d\mathbf{x}).
		\end{aligned}
		\]
		
		(2) In the case $0<s\le 1/2$, the same compactness result also holds in case $\{u_{i}\}$ is a sequence of minimizing $s$-harmonic	maps.	Moreover, in this case the limit $u$ is also a minimizing $s$-harmonic map.
	\end{theorem}
This compactness result implies  that if  $u \in \widehat{H}^s\left(\Omega ; \mathbb{S}^{d-1}\right)$ is a stationary $s$-harmonic map for $s\neq 1/2$, then for every $x\in \Om$  and every sequence $r_k\to 0$, there exists a subsequence  $r^{\prime}_k\to 0$ and a map $\var:\R^n\to \mathbb{S}^{d-1}$ which is 0-homogeneous at the origin (see Definition \ref{def: classical symmetry} below) such that
	\[
	\begin{aligned} & u_{x,r_k^\prime}\to \var\qquad\,\,\text{strongly in }\widehat{H}^{s}(D_R, \mathbb{R}^d), \text{  and  }\\
		& u_{x,r_k^\prime}^{e}\to \var^{e}\qquad\text{strongly in }H^{1}(B_R^+;\mathbb{R}^d,|z|^{a}d\mathbf{x})
	\end{aligned}
	\]
	for all $R>0$, where $u_{x,r}(y)= u(x+ry)$. For details, see \cite[Section 7.2]{Millot-Pegon-Schikorra-2021-ARMA}.
	
	\begin{definition} Let  $s\neq 1/2$ and $u \in \widehat{H}^s\left(\Omega ; \mathbb{R}^d\right)$ be a stationary $s$-harmonic map. The above deduced $0$-homogeneous map $\var$ is called a tangent map of $u$ at $x\in \Om$. \end{definition}	
	
As that of \cite{Cheeger-Naber-2013-CPAM, Naber-V-2017}, the symmetry of tangent maps motivates the  notion of symmetric functions so as to stratify singular set of a $s$-harmonic map.
	
	\begin{definition}[$k$-symmetry]\label{def: classical symmetry} Given a measurable map $\var: \mathbb{R}^n \rightarrow \mathbb{R}$. We say that
		
		(1)  $\var$ is  $0$-homogeneous or  $0$-symmetric with respect to point $p \in \mathbb{R}^n$ if $\var(p+\lambda v)=$ $\var(p+v)$ for all $\lambda>0$ and $v \in \mathbb{R}^n$.
		
		(2)  $\var$ is $k$-symmetric if $\var$ is $0$-symmetric with respect to the origin, and $\var$ is translation invariant with respect to a $k$-dimensional subspace $V \subset \mathbb{R}^n$, i.e.,
		$$
		\var(x+v)=\var(x) \quad \text { for all } x \in \mathbb{R}^n, v \in V.
		$$
	\end{definition}
	
	Then, for any $k\in\{0,1,\cdots,n\}$, we can  define the set
	\begin{equation}\label{eq: classical stratification}
	\Sigma^{k}(u)=\{x\in\Om:\text{no tangent map of }u\text{ is }(k+1)\text{-symmetric at }x\}.
	\end{equation}
such that
	\[\Sigma^{0}(u) \subset \Sigma^{1}(u)\subset \cdots\subset \Sigma^{n-1}(u)\subset \Sigma^{n}(u)=\Om. \]
	Now we have the simple observation.
	\begin{remark}
	 Suppose $s\in (0,1)\backslash \{1/2\}$ and $u \in \widehat{H}^s\left(\Omega ; \mathbb{S}^{d-1}\right)$ is a stationary $s$-harmonic map. Then
	$$\sing(u)= \Sigma^{n-1}(u).$$
	As a result,
	 \[\Sigma^{0}(u) \subset \Sigma^{1}(u)\subset \cdots\subset \Sigma^{n-1}(u)=\sing(u). \]
	\end{remark}
This is the so-called classical stratification of $\sing(u)$.
	To see this, suppose that $u$ has a constant tangent map at $x$,
	that is,  $x\in \Om\backslash \Sigma^{n-1}(u)$. Theorem \ref{thm: compactness} implies that $\boldsymbol{\Theta}_s(u^e, \mathbf{x}, r)\to 0$ as $r\to 0$,  which in turn implies that $u$ is smooth in a neighborhood of $x$ by the $\ep$-regularity theorem. Hence $\sing(u)\subset \Sigma^{n-1}(u)$.  On the other hand, if $u$ is smooth in a neighborhood of $x$, then surely there is a unique constant tangent map at $x$, and thus $ \Sigma^{n-1}(u)\subset\sing(u)$. This proves the remark.

\section{Quantitative stratification of singular set}

\subsection{The setting of quantitative symmetry }
Now we recall the notion of quantitative symmetry introduced already in our previous paper \cite{He-X-Z-24}. Throughout, we will always assume $s\in (0,1)\backslash \{1/2\}$ and fix a constant $$1<\ga_0<\min\{1/(1-s), 2\}.$$
\begin{definition}[Quantitative symmetry]\label{quantitative symmetry}  Given a map $u \in \widehat{H}^{s}(\Omega), \ep>0$ and nonnegative integer $k$, we say that $u$ is $(k, \ep)$-symmetric on $D_r(x) \subset \Omega$, if there exists a  $k$-symmetric function $h\in  \widehat{H}^{s}(D_{r}, \R^d)$ such that
	$$	
	\medint_{B_1^{+}}\left|u^e_{x,r}(\mathbf{y})
	-h^e_{0,r}(\mathbf{y}) \right|^{\gamma_0} d \mathbf{y}   =	\medint_{B_r^{+}(\mathbf{x})}\left|u^e(\mathbf{y})
	-h^e(\mathbf{y}-\mathbf{x})\right|^{\gamma_0} d \mathbf{y} \leq \ep,
	$$
	where $\mathbf{x}=(x, 0)$.\end{definition}
Note that the above quantitative symmetry is defined via the upper space extension, not as that of harmonic maps.

\begin{remark}Equivalently, $u$ is $(k,\ep)$-symmetric on $D_{r}(x)$ if and only
	if the scaled map $u_{x,r}(y)=u(x+ry)$ is $(k,\ep)$-symmetric on
	$D_{1}$. The choice of $\ga_0$ is used to make sure that $H^1(B_r^+, z^ad\mathbf{x})$ is compactly embedded into  $L^{\ga_0}(B_r^+)$ for all $s\in (0,1)$ (see \eqref{eq: comp-embedding}).
\end{remark}
 The following proposition is an important property of this quantitative symmetry.

\begin{proposition}\label{rmk: compactness-of-symmetry}
	Suppose $\{u_i\}\subset \widehat{H}^{s}(D_2)$ converges weakly to some function $u\in \widehat{H}^{s}(D_2)$. If $u_i$ is $(k,\ep_i)$-symmetric in $D_1$ for some $\ep_i\to 0$, then $u$ is $k$-symmetric in  $D_1$ as well.
\end{proposition}

The proof can be referenced in \cite[Remark 3.2]{He-X-Z-24}, and so omitted.
Now we recall  the definition of quantitative symmetry introduced in \cite[Remark 3.2]{He-X-Z-24}.

\begin{definition}[Quantitative stratification]\label{ddd} For any
	map  $u\in\widehat{H}^{s}(\Omega,N)$, $r,\ep>0$ and $k\in\{0,1,\cdots,n\}$,
	we define the $k$-th quantitative singular stratum $S_{\ep,r}^{k}(u)$
	by
	\[
	S_{\ep,r}^{k}(u)\equiv\big\{ x\in\Omega:u\text{ is not }(k+1,\ep)\text{-symmetric on }D_{s}(x)\text{ for any }r\le s<1\big\}.
	\]
	Furthermore, we set
	\[
	S_{\ep}^{k}(u):=\bigcap_{r>0}S_{\ep,r}^{k}(u)\quad\text{ and }\quad S^{k}(u)=\bigcup_{\ep>0}S_{\epsilon}^{k}(u).
	\]
\end{definition}

It follows that
\begin{equation}\label{Skepr}
k'\leq k\text{ or }\ep'\geq\ep\text{ or } r'\leq r\ \Longrightarrow \ S^{k'}_{\ep',r'}(u)\subseteq S^k_{\ep,r}(u).
\end{equation}
In particular,
\[
S^0(u)\subset S^1(u)\subset\cdots\subset S^n(u)=\Omega.
\]
The following lemma shows that $S^k(u)$ is well adapted to stationary $s$-harmonic maps and was proved by \cite[Remark 3.4]{He-X-Z-24}.

\begin{lemma}\label{equi}   Let $s\in (0,1)\backslash \{1/2\}$ and  suppose $u\in \widehat{H}^{s}(\Om,\mathbb{S}^{d-1})$ is a stationary $s$-harmonic map. Then
	\[
	S^{k}(u)=\Sigma^k(u)
	\]
	where $\Sigma^k(u)$ is defined as in \eqref{eq: classical stratification}.
	Consequently, we have
	\[
	S^{0}(u)\subset S^{1}(u)\subset\cdots\subset S^{n-1}(u)={\rm sing}(u).
	\]
\end{lemma}

%

 Observe that if  $u$ is a stationary  $s$-harmonic map and $\boldsymbol{\Theta}_s(u^e,\mathbf{x},1)=\boldsymbol{\Theta}_s(u^e,\mathbf{x},1/2)$, then by the monotonicity formula \eqref{1.8} and the unique continuation lemma \ref{lem: unique continuation property}, $u^e$ must be $0$-symmetric with respect to $\textbf{x}$ in $\R^{n+1}_+$. The following proposition gives a sufficient condition for $(0,\ep)$-symmetry of a stationary  $s$-harmonic map.
\begin{proposition}\label{pro1}	 Let $s\in (0,1)\backslash \{1/2\}$. For any $\ep>0,$ there exists $\delta_1=\delta_1(n,N,\Lambda,s,\ep)$ such that, if  $u\in \widehat{H}^{s}_\La(D_{10},N)$ (see \eqref{eq: bounded energy space} for definition) is a stationary weakly $s$-harmonic map with
	$$\boldsymbol{\Theta}_s(u^e,\mathbf{x},r)
-\boldsymbol{\Theta}_s(u^e,\mathbf{x},r/2)<\delta_1$$
	for some $x\in D_1$ and $0<r<1$, then $u$ is $(0,\ep)$-symmetric on $D_r(x)$.
\end{proposition}
\begin{proof}
	Consider by contradiction a sequence of  stationary $s$-harmonic  maps $u_i\in \widehat{H}^{s}_\La(D_{10},N)$ which satisfy
	$$\boldsymbol{\Theta}_s(u^e,\mathbf{x}_i,r_i)
-\boldsymbol{\Theta}_s(u^e,\mathbf{x}_i,r_i/2)<i^{-1},$$  but are not $(0,\ep)$-symmetric on $D_{r_i}(x_i)$.
	Let $\bar{u}_i(y)=u_i(x_i+r_iy)$. Then we have $\bar{u}^e_i(\textbf{y})=u^e_i(\textbf{x}_i+r_i\textbf{y})$ and
	$$\boldsymbol{\Theta}_s(\bar{u}_i^e,\mathbf{0},\tau) =\boldsymbol{\Theta}_s(u_i^e,\mathbf{x}_i,\tau r_i),\qquad \forall\,\tau>0.$$
	By the monotonicity formula \eqref{1.8}, $\{\bar{u}^e_i\}_{i\geq 1}$ is uniformly bounded in $H^1\left(B_1^+,|z|^a \mathrm{d} \mathbf{x}\right)$. So up to  a subsequence, we may  assume that $\bar{u}^e_i\rightharpoonup v^e$ weakly in $H^1\left(B_1^+,|z|^a \mathrm{d} \mathbf{x}\right)$. Using the monotonicity formula \eqref{1.8} again, for all $i\in \N$, we have
	\[
	\boldsymbol{\de}_s\int_{B^+_1\backslash B^+_{1/2}}z^a\frac{|\mathbf{x}\cdot\nabla \bar{u}^e_i|^2}{|\mathbf{x}|^{n+2-2s}}d\mathbf{x} \le \boldsymbol{\Theta}_s(\bar{u}_i^e,\mathbf{0},1) -\boldsymbol{\Theta}_s(\bar{u}_i^e,\mathbf{0},1/2)\le 1/i.
	\]
Sending $i\to \infty$, we deduce that $v^e$ is radially invariant on $B^+_1\backslash B^+_{1/2}$. Again, this implies that $v^e$ is 0-homogeneous on $\R^{n+1}$ by the unique continuation property. Finally, since $\bar{u}^e_i\rightarrow v^e$ strongly in $L^{\gamma_0}(B_1^+)$, it follows that
	\[
	\medint_{B^+_{r_i}(\textbf{x}_i)} |u^e_i-v^e((\cdot-\textbf{x}_i)/r_i)|^{\gamma_0}
=\medint_{B^+_{1}}|\bar{u}^e_i-v^e|^{\gamma_0}\to 0\qquad \text{as }i\to \infty,
	\]
	contradicting with the assumption that $u_i$ is not  $(0,\ep)$-symmetric on $D_{r_i}(x_i)$ for all $i$.
\end{proof}

\subsection{Properties of quantitative stratum} Before discussing stratified singular stratum, we first recall the notion of  quantitative frame for affine subspaces as in \cite{Naber-V-2018}.
\begin{definition}
	Let $\{y_i\}_{i=0}^k$ and $\rho>0$. We say that $\{y_i\}_{i=0}^k$ is a $\rho$-frame, or we say that these points $\rho$-effectively span a $k$ dimensional affine subspace, if for all $i=1,\cdots,k$,
	$$\dist(y_i,y_0+\Span\{y_1-y_0,\cdots,y_{i-1}-y_0\})\geq2\rho.$$
	We say that a set $F$ spans  $\rho$-effectively  a $k$ dimensional affine subspace, if there exist points $\{y_i\}_{i=0}^k\subset F$ which spans $\rho$-effectively  a $k$ dimensional affine subspace.
\end{definition}

\begin{remark}\label{rmk: quantitative frame} (\cite[Definition 28]{Naber-V-2018})
	The advantage of the above introduced quantitative frame is twofold:
	
	(1)  If $\{y_i\}_{i=0}^k$ $\rho$-effectively spans a $k$-dimensional affine subspace, then for every point 	$x\in  y_0+{\rm span}\{y_1-y_0,\cdots,y_k-y_0\}$,
	there exists a unique set of numbers $\{\alpha_i\}_{i=1}^k$ such that
	$$x=y_0+\sum_{i=1}^k\alpha_i(y_i-y_0) \quad \text{with}\quad  |\alpha_i|\leq C(n,\rho)|x-y_0|.$$
	
	(2) Quantitative frame is stable under limiting process: if, for each $j\ge 1$, the set $\{y_i^j\}_{i=0}^k$ $\rho$-effectively spans a $k$ dimensional affine subspace and $y_i^j\to y_i$ as $j\to \infty$, then $\{y_i\}_{i=0}^k$ also $\rho$-effectively spans a $k$ dimensional affine subspace.
\end{remark}

The next proposition shows that $\rho$-frame can extend the above rigidity result of Proposition \ref{pro1} to higher dimension, which is also called quantitative cone splitting in the literature.

\begin{proposition} (Quantitative cone splitting). For any $\epsilon, \rho>0$, there exists $\delta=$ $\delta(n, d, \Lambda, \ep, \rho)>0$ satisfying the following property. Let $u \in \widehat{H}_{\Lambda}^s\left(D_{10}, \mathbb{S}^{d-1}\right)$ be a stationary $s$-harmonic map, and there exist points $\left\{x_j\right\}_{j=0}^k \subset D_1$ which satisfies

(i) $\boldsymbol{\Theta}_s(u^e,\mathbf{x}_j,1)
-\boldsymbol{\Theta}_s(u^e,\mathbf{x}_j,1/2)<\delta$ for all $j=0, \cdots, k$,

(ii) The set $\left\{x_j\right\}_{j=0}^k $ $\rho$-effectively spans a $k$-dimensional affine plane $V$.

Then $u$ is $(k, \varepsilon)$-symmetric on $D_1$.

\end{proposition}

The proof of this proposition is similar to that of Proposition \ref{pro1} and thus omitted. Using the notion of quantitative frame, the next proposition shows that for a sufficiently pinched singular stratum of a stationary $s$-harmonic mapping, if it is of high dimension in essential, then it  satisfies a  one-side Reifenberg approximating property.
\begin{proposition}\label{pro32} Let $s\in (0,1)\backslash\{1/2\}$.
	For any $1>\ep,\rho>0$, there exists a positive small constant $\delta_2=\delta_2(n,N,\Lambda,s,\ep,\rho)$ such that the following holds. For any stationary  $s$-harmonic map $u\in \widehat{H}^{s}_{\Lambda}(D_{10},N)$, suppose  that the collection of points
	$$F_{\delta_2}:=\left\{y\in D_2:  \boldsymbol{\Theta}_s(u^e,\mathbf{y},\rho) >E-\delta_2\right\}$$
	$\rho$-effectively spans a $k$-dimensional affine plane $V\subset\R^n$, where
	$$E:=\sup_{y\in  D_2}\boldsymbol{\Theta}_s(u^e,\mathbf{y},4).$$
	Then
	$$S_{\ep,\delta_2}^k\cap D_1\subset D_{2\rho}(V).$$

\end{proposition}
To prove this proposition, we first give a technical lemma. The notation below turns out to be quite useful: For any $k$ dimensional subspace $L$, we denote
$$|L\cdot\nabla f|^2=\sum\limits_{i=1}^{k}|\nabla_{e_i}f|^2,$$
where $\{e_i\}_{i=1}^{k}$ is an arbitrary orthonormal basis of $L$. This quantity measures how far away the function $f$  is translation invariant along the given subspace $L$.
\begin{lemma}\label{lem: technical lemma} Let $s\in (0,1)\backslash\{1/2\}$ and let $\rho,\ep>0$ be fixed. Then there exists $\tau=\tau(n,\La,\rho,\ep,s)>0$ satisfying the following property. Suppose $u\in \widehat{H}^{s}_{\Lambda}(D_{10},N)$ is  a stationary  $s$-harmonic map and $$\int_{B^+_1}|P\cdot \na u^e|^2 z^a d\mathbf{x}<\tau,$$
for some $k+1$ dimensional plane $P$ in $\R^n$, then $S^k_{\ep, \bar{\tau}}\cap D_{1/2}=\emptyset$, where $\bar{\tau}=\tau^{1/2(n-2s)}$.
\end{lemma}
 The proof is totally the same as that of Lemma 32 of Naber-Valtorta \cite{Naber-V-2018} and thus omitted.  Now we prove Proposition \ref{pro32}.

\begin{proof}
	By assumption, we assume that  $V$	is  the affine subspace that is spanned by a $\rho$-frame  $\{y_{i}\}_{i=0}^{k}\subset F_{\delta_2}$, i.e.,
	\[
	V=y_{0}+{\rm span}\{y_{i}-y_{0}\}_{i=1}^{k}.
	\]
	Let $\mathbf{x}_{0}=(x_0,0),x_0\in D_{1}\backslash D_{2\rho}(V)$ and let $\delta_2>0$ to be determined later. Then, for $0<r<\rho$ there are
	$$B^+_{r}(\mathbf{x}_0)\subset B^+_{4}(\mathbf{y}_i)\backslash B^+_{\rho}(\mathbf{y}_i),\qquad \forall\,0\le i\le k.$$

We claim that there exist constants $r=r(n,\La,s,\rho)>0$ sufficiently small and $C_1=C_1(n,\Lambda,s,\rho)>0$ such that
\begin{equation}\label{eq: small translation}
 	r^{2s-n}\int_{B_r^+(\mathbf{x}_0)}z^a|L\cdot\na u^e(\mathbf{x})|^2d\mathbf{x}< \tau
\end{equation}
	for some $k+1$ dimensional subspaces $L\subset \mathbb{R}^n$, where $\tau$ is defined as in Lemma \ref{lem: technical lemma}.
	
	To prove this claim, first note that by the definition of $F_{\delta_2}$, there holds $$\boldsymbol{\Theta}_s(u^e,\mathbf{y}_i,4)
-\boldsymbol{\Theta}_s(u^e,\mathbf{y}_i,\rho)<\de_{2},\quad \forall\,0\le i\le k.$$
	Then, by the monotonicity formula \eqref{1.8} and the definition of $F_{\de_2}$, we have
	\[
	\int_{B^+_{r}(\mathbf{x}_0)}z^a|(\textbf{x}-\mathbf{y}_i)\cdot\na u^e(\mathbf{x})|^{2}d\textbf{x}\le C\int_{B^+_{4}(\mathbf{y}_i)\backslash B^+_{\rho}(\mathbf{y}_i)}z^a\frac{|(\textbf{x}-\mathbf{y}_i)\cdot\na u^e(\mathbf{x})|^{2}}{|\mathbf{x}-\mathbf{y}_i|^{n+2-2s}}d\textbf{x}\le C_{n,\rho,s}\de_{2}.
	\]
	Consequently, by an elementary triangle inequality we deduce
	\[
	\int_{B^+_{r}(\mathbf{x}_0)}z^a|(\mathbf{y}_i-\mathbf{y}_0)\cdot\na u^e(\mathbf{x})|^{2}d\textbf{x} \le C_{n,\rho,s}\de_{2},\quad \forall\,1\le i\le k.
	\]
	Since $\{y_{j}\}_{j=0}^{k}$ is $\rho$-independent, we conclude that
	\begin{equation}
		\int_{B^+_{r}(\mathbf{x}_0)}z^a|V\cdot\na u^e(\mathbf{x})|^{2}d\textbf{x}\le C(n,\rho,s)\de_{2}\label{eq: k-small}
	\end{equation}
	where $V={\rm span}\{ {y}_i-y_{0}\}_{i=1}^{k}\subset \R^n$.
	
	On the other hand, for each $\mathbf{x}\in B_{r}^+(\mathbf{x}_0)\subset B_{2}^+\backslash B_{\rho}^+(V)$,
	let
	\[
	\pi_{V}(\mathbf{x})=y_{0}+\sum_{i=1}^{k}\al_{i}(\mathbf{x})(y_{i}-y_{0})
	\]
	be the orthogonal projection of $\mathbf{x}$ on $V$. We have $| \mathbf{x}-\pi_{V}( \mathbf{x})|\ge\rho$,
	$|\al_{i}(\mathbf{x})|\le C(n,\rho)$ (see Remark \ref{rmk: quantitative frame}) and
	\[
	\begin{aligned}
	\int_{B^+_{r}(\mathbf{x}_0)}z^a|(\textbf{x}-\pi_{V}(\textbf{x}))\cdot\na u^e(\mathbf{x})|^{2}d\textbf{x} &\le C\int_{B_r^+(\mathbf{x}_0)}z^a|(\textbf{x}-\mathbf{y}_{0})\cdot\na u^e(\mathbf{x})|^{2}d\textbf{x}\\
		&\quad  + C(n, \rho)\sum_{i=1}^{k}\int_{B_r^+(\mathbf{x}_0)}z^a|(\mathbf{y}_{i}-\mathbf{y}_{0})
\cdot\na u^e(\mathbf{x})|^{2}d\textbf{x}\\
		&\le C(n,\rho,s)\de_{2}.
	\end{aligned}
	\]
	Thus, by setting $\Large{h(\mathbf{x})=\frac{\mathbf{x}-\pi_{V}(\mathbf{x})}{|\mathbf{x}-\pi_{V}
(\mathbf{x})|}}$, it follows	that
	\[
	\begin{aligned}\int_{B^+_{r}(\mathbf{x}_0)}z^a|h(\mathbf{x}_0)\cdot\na u^e(\mathbf{x})|^{2}d\textbf{x} & \le C\int_{B^+_{r}(\mathbf{x}_0)}z^a|h(\textbf{x})\cdot\na u^e(\mathbf{x})|^{2}\\ &\quad+C\int_{B^+_{r}(\mathbf{x}_0)}z^a|(h(\textbf{x})-h(\mathbf{x}_0))\cdot\na u^e(\mathbf{x})|^{2}d\textbf{x}\\
		& \le C\de_{2}+Cr^{2}\int_{B^+_{r}(\mathbf{x}_0)}z^a|\na u^e(\mathbf{x})|^{2}d\textbf{x}\\
		& \le C\de_{2}+C(n,\rho,s,\La)r^{n+2-2s}.
	\end{aligned}
	\]
Put $L:=V\oplus \R h(\mathbf{x}_0)\subset\R^n$ be a $k+1$ dimensional subspace. Then, combining the above estimate together with \eqref{eq: k-small},  we find that
\[
	r^{2s-n} \int_{B^+_{r}(\mathbf{x}_0)}z^a|L\cdot\na u^e(\mathbf{x})|^{2}d\textbf{x}\le C(n,\Lambda,s,\rho)\de_{2}r^{2s-n}+Cr^{2}.
	\]
Now we choose $r=r(n,\La,\rho,s)\ll\rho$ such that such that
$Cr^2<\tau/2$,
and then choose $\de_2$ such that
$C(n,\Lambda,s,\rho)\delta_2 r^{2s-n}<\tau/2$.
The estimate \eqref{eq: small translation} and thus the claim follows.
	
	Finally, applying Lemma \ref{lem: technical lemma} to  \eqref{eq: small translation}  we find that $x_0\not\in S^k_{\ep,\bar{\tau}r}$. Easy to note that $\de_2\ll \bar{\tau}r$. Hence using the conclusion $S^k_{\ep,\de_2}\subset S^k_{\ep,\bar{\tau}r}$ we conclude that $x_0\not\in S^k_{\ep,\delta_2}(u)$. The proof is complete.
\end{proof}
The next result shows that $\boldsymbol{\Theta}_s(u^e,\cdot,\rho)$ remains almost constant on all pinched points.
\begin{proposition}\label{5.5}
	For any $0<\rho,\eta<1$,  there exists $\de_3=\de_3(n,d,\Lambda,s,\rho,\eta)>0$  satisfying the following property. Let  $u\in \widehat{H}^{s}_{\Lambda}(D_{10},\mathbb{S}^{d-1})$ be a stationary  $s$-harmonic map, if the set
	$$
	F_{\delta_3}:=\{y\in  D_2:\boldsymbol{\Theta}_s(u^e,\mathbf{y},\rho)>E-\delta_3\}
	$$
	$\rho$-effectively spans a $k$ dimensional affine subspace $L\subset\mathbb{R}^n$, then
	$$\boldsymbol{\Theta}_s(u^e,\mathbf{x},\rho) \ge E-\eta \quad \text{for all }  x\in L\cap D_1,$$
	where $E$ is defined in Proposition \ref{pro32}.
\end{proposition}
\begin{proof}
	Suppose, for some $\rho,\eta>0$, there exists a sequence of stationary $s$-harmonic maps $\{u_i\}_{i\geq 1}\subset \widehat{H}^{s}_{\Lambda}(D_{10},\mathbb{S}^{d-1})$ with
	$$
	\sup_{y\in  D_2}\boldsymbol{\Theta}_s(u^e_i,\mathbf{y},4) \leq E_i.
	$$
	For each $i\ge 1$, the set $F_{i}:=\{y\in D_2:\boldsymbol{\Theta}_s(u^e_i,\mathbf{y},\rho) >E_i-i^{-1}\}$ contains a subset $\{y^i_j\}_{j=0}^k$ spanning $\rho$-effectively  a $k$ dimensional affine subspace $L_i\subset\mathbb{R}^n$, and there exists $x_i\in L_i\cap D_1$ such that
	\begin{align}\label{fanz}
		\boldsymbol{\Theta}_s(u^e_i,\mathbf{x}_i,\rho) \leq E_i-\eta.
	\end{align}
	
	The first consequence of the assumption is that
	\begin{equation}\label{eq:lemma 312} \boldsymbol{\Theta}_s(u^e_i,\mathbf{y}_j^i,4)
-\boldsymbol{\Theta}_s(u^e_i,\mathbf{y}_j^i,\rho) <1/i,\quad \forall\,i\ge 1\text{ and } 0\le j\le k.
	\end{equation}
	Without loss of generality, we  further assume that, for each $0\le j\le k$,
	\[y^i_j\to y_j,\,\, x_i\to x_0\,\,\text{and}\,\, E_i\to E_0 \quad \text{as }i\to \infty\]
	and that the sequence  $L_i$ of $k$ dimensional affine subspaces converges to a $k$ dimensional affine subspace $L$ passing through  $x_0$ and $\{y_j\}_{0}^k$.
	
	By the compactness Theorem \ref{thm: compactness}, there exist a stationary $s$-harmonic map $v$ such that
	\[
	\begin{cases}
	u_i\to v \quad &\text{in} \,\, \widehat{H}^s(D_4,\mathbb{R}^d),\\
	u^e_i\to v^e \quad  &\text{in} \,\,H^1(B_4^+, |z|^ad\mathbf{x}).
	\end{cases}
	\]
	Using the monotonicity formula \eqref{1.8} and the stationarity of $v$, we find that for any fixed $\mathbf{y}=\{y,0\}$, the function
	$$
	r\mapsto \boldsymbol{\Theta}_s(v^e,\mathbf{y},r)
	$$
	is monotonically nondecreasing. Thus sending $i\to \infty$ in \eqref{fanz} and \eqref{eq:lemma 312}, we obtain
	\[
	\boldsymbol{\Theta}_s(v^e,\mathbf{y}_j,4)
=\boldsymbol{\Theta}_s(v^e,\mathbf{y}_j,\rho)  =E_0,\quad \text{for all }\, 0\le j\le k
	\]
	and \begin{equation}\label{eq: loss of energy}
		\boldsymbol{\Theta}_s(v^e,\mathbf{x}_0,\rho) \le E_0-\eta.
	\end{equation}
	Moreover, by the unique continuation property and classical cone splitting principle, we know that  $v^e$ is  translation invariant along $L$. Hence
	$$
	\boldsymbol{\Theta}_s(v^e,\mathbf{x}_0,\rho) \equiv \boldsymbol{\Theta}_s(v^e,\mathbf{y}_j,\rho)=E_0,
	$$
	which clearly contradicts with \eqref{eq: loss of energy}. The proof is complete.
\end{proof}

The following technical result shows that the almost symmetry is preserved under certain pinching condition.
\begin{proposition}\label{Non-symmetry} Let $s\in (0,1)\backslash\{1/2\}$. For any $\ep>0$, there exist $\ep_{\ast},\de_4>0$,
depending only on $n,\Lambda,d,s,\ep$, satisfying the following property.
Let $u\in\widehat{H}_{\Lambda}^{s}(D_{10},\mathbb{S}^{d-1})$ be a
stationary $s$-harmonic map and $0<r<1$. Let $S\subset S_{\ep,\de_{4}r}^{k}(u)\cap D_{1}$
and $x\in D_{2}$. If $D_{t}(x)\cap S\neq\emptyset$ holds for some
$1\ge t\ge r$, then $x\in S_{\ep_{\ast},t}^{k}(u)$.\end{proposition}
\begin{proof}
Suppose, on the contrary, there exists $\ep>0$ such that, for any
$\ep_{i},\de_{i}\to0$, there are stationary $s$-harmonic maps $u_{i}\in\widehat{H}_{\Lambda}^{s}(D_{10},\mathbb{S}^{d-1})$,
$0<r_{i}<t_{i}\le\rho_{i}<1$, $S_{i}\subset S_{\ep,\de_{i}r_{i}}^{k}(u_{i})\cap D_{1}$
and $x_{i}\in D_{2}$ satisfying $p_{i}\in D_{t_{i}}(x_{i})\cap S_{i}\neq\emptyset$
and $u_{i}$ is $(k+1,\ep_{i})$-symmetric on $D_{\rho_{i}}(x_{i})$.

Put $\bar{u}_{i}(x)=u_{i}(x_{i}+\rho_{i}x)$. It follows from the monotonicity formula that $\bar{u}_{i}$
is uniformly bounded on $D_{5}$ and $(k+1,\ep_{i})$-symmetric
on $D_{1}(0)$. Thus, by compactness theorems, we may assume, up to
a subsequence, that $\bar{u}_{i}\to v$ strongly in $\widehat{H}^{s}(D_{4},\mathbb{S}^{d-1})$
and $\bar{u}_{i}^{e}\to v^{e}$ strongly in $H^{1}(B_{4}^{+},z^{a}d\mathrm{x})$.
Moreover, $v$ is $k+1$ symmetric on $D_{1}$. Thus there exists a $k+1$ plane $P\subset\R^{n}$
such that
\[
P\cdot\na v^{e}=0 \qquad\text{in }B^+_1.
\]
Thus for any vector $w\in P$, we thus have $w\cdot \na v^e=0$ in $B^+_1$. Note that $\tilde{v}=w\cdot \na v^e$ solves the linear equation $\mathrm{div}(z^a \na \tilde{v})=0$. Hence the unique continuation principle implies that $\tilde{v}\equiv 0$ in $B^+_4$. In particular, we have $P\cdot\na v^{e}=0 \qquad\text{in }B^+_4.$

Next, let $q_{i}=(p_{i}-x_{i})/\rho_{i}$. Then $|q_{i}|\le1$ and
\begin{equation}
q_{i}\in\frac{S_{i}-x_{i}}{\rho_{i}}\cap D_{\frac{t_{i}}{\rho_{i}}}\subset\frac{1}{\rho_{i}}\left(S_{\ep,\de_{i}r_{i}}^{k}(u_{i})\cap D_{1}-x_{i}\right)\cap D_{\frac{t_{i}}{\rho_{i}}}\subset S_{\ep,\frac{\de_{i}r_{i}}{\rho_{i}}}^{k}(\bar{u}_{i})\cap D_{1}.\label{eq: XX}
\end{equation}
We may assume up to a subsequence that $q_{i}\to q\in\bar{D}_{1}$.
Fix a  small constant $h>0$. Then the strong convergence
property implies that
\[
h^{2s-n}\int_{B_{h}^{+}(\mathrm{q}_{i})}z^a|P\cdot\na\bar{u}_{i}^{e}|^{2}\to0\qquad\text{as }i\to\wq.
\]
 As a result, Lemma \ref{lem: technical lemma} implies that $q_{i}\not\in S_{\ep,c_{h}}^{k}(\bar{u}_i)$
for  $c_{h}=\bar{\tau}h>0$ depending only on a scaling of $h$, independent of $i$. However,
note that $\de_{i}r_{i}/\rho_{i}\le\de_{i}<c_{h}$ as $i\to\wq$.
We reach a contradiction with (\ref{eq: XX}). The proof is complete.
\end{proof}

\section{Reifenberg theorems and estimates of Jones' number}

Reifenberg theorem plays a key role in the proof of rectifiability theorems, in which a key quantity is Jones' number. Following the line of Naber-Valtorta \cite{Naber-V-2018}, we first recall Jones'  number  $\beta_2$, which quantifies how close the support of a measure $\mu$ is to a $k$ dimensional affine subspace, and then recall the Reifenberg type results obtained by Naber-Valtorta~\cite{Naber-V-2017}. Then we will extend the estimate of Jones' number of \cite{Naber-V-2017} to our case.

\begin{definition}\label{aaa}
	Let $\mu$ be a nonnegative Radon measure defined on $D_3$. Fix  $k\in\N$, the $k$ dimensional Jones' number  $\beta_\mu^k$ is defined as
	\begin{align*}
		\beta^k_{\mu}(x,r)=\bigg(\inf_{V\su\R^n}\int_{D_r(x)}\frac{d^2(y,V)}{r^2}\frac{d\mu(y)}{r^k}\bigg)^{\frac{1}{2}},\quad \text{when }\,D_r(x)\subset D_3,
	\end{align*}
	where the infimum is taken over  all $k$ dimensional  affine subspaces $V$.
\end{definition}

It turns out that  Jones' number can be calculated by making use of probability theory, see below.
We now state two versions of quantitative Reifenberg theorems from \cite{Naber-V-2017}. The reader can immediately note that Jones' number stays in the core of the below two theorems.
\begin{theorem}[{Discrete-Reifenberg, \cite[Theorem 3.4]{Naber-V-2017}}]\label{rre}
	There exist two constants $\de_5=\de_5(n)$ and $C_R(n)$ such that the following property holds. Let $\{D_{r_x}(x)\}_{x\in \mathcal{C}}\su D_2\su\R^n$ be a family of pairwise disjoint balls with centers in $\mathcal{C}\subset D_1$ and let $\mu\equiv\sum_{x\in \mathcal{C}}\omega_kr_x^k\de_x$ be the associated measure. If for every ball $x\su D_1$ and $0<r<1$ there holds
	\begin{align}\label{eq: multiscale appro-1}
		\int_{D_r(x)}\left(\int_0^r\beta^k_{\mu}(y,s)^2\frac{ds}{s}\right)d\mu(y)<\de_5^2r^k,
	\end{align}
	then we have the uniform estimate
	\begin{align*}
		\sum_{x\in \mathcal{C}}r^k_x< C_R(n).
	\end{align*}
\end{theorem}
The second Reifenberg type theorem reads as follows.

\begin{theorem} [{Rectifiable-Reifenberg, \cite[Theorem 3.3]{Naber-V-2017}}]\label{rree}
	There exist constants $\de_6=\de_6(n)$ and $C=C(n)$ satisfying the following property. Assume that $S\subset D_2\su\R^n$ is $\HH^k$-measurable, and for each $D_r(x)\subset D_2$ there holds
	\begin{align}\label{eq: multiscale appro-2}
		\int_{S\cap D_r(x)}\left(\int_0^r\beta^k_{\HH^k|_S}(y,s)^2\frac{ds}{s}\right)d\HH^k(y)<\de_6^2r^k.
	\end{align}
	Then $S\cap D_1$ is  $k$-rectifiable, and   $\HH^k(S\cap D_r(x))\leq C r^k$ for each $x\in S\cap D_1$.
\end{theorem}

The two conditions \eqref{eq: multiscale appro-1} and \eqref{eq: multiscale appro-2} are usually called multiscale approximation conditions. We point out that Mi\'skiewicz \cite{Miskiewicz-18-Finn} obtained some  improvement for the above two Reifenberg theorems.

To apply the above two theorems, we need to estimate Jones' number $\beta_2$. Following the approach of   \cite[Theorem 7.1]{Naber-V-2017} for harmonic maps,  we need to establish an $L^2$ subspace approximation theorem. To this end, for $x\in D_1, \mathbf{x}:=(x,0) $ and $r>0,$ we define
\begin{equation}\label{W_rx}
	W_r({x}):= W_{r,8r}({x})=\int_{B^+_{8r}({\mathbf{x}})\ba B^+_r({\mathbf{x}})}z^a\frac{|(\mathbf{y}-\mathbf{x})\cdot\na u^e(\mathbf{y})|^2}{|\mathbf{y}-\mathbf{x}|^{n+2-2s}}d\mathbf{y}.
\end{equation}

\begin{theorem}\label{ji} Let $s\in (0,1)\backslash \{1/2\}$.
	For any $\ep>0$, there exist $C(n,\Lambda,s,\ep)>0$ and $\de_7(n,\Lambda,s,\ep)>0$ such that, for any  stationary $s$-harmonic map $u\in \widehat{H}^{s}_\La(D_{10},\mathbb{S}^{d-1})$ and $0<r\leq1$,  $x\in D_1$, if $u$ is  $(0,\de_7)$-symmetric  but not $(k+1,\ep)$-symmetric on $D_{8r}(x)$, then for any nonnegative finite measure $\mu$ on $D_r(x)$,  we have
	\begin{align}\label{po}
		\beta^k_{\mu}({x},r)^2\leq Cr^{-k}\int_{D_r(x)}W_r({y})d\mu({y}).
	\end{align}
\end{theorem}
In order to prove this theorem, we can assume that $x=0,r=1$ and $\mu$ is  a probability measure supported on  $D_1$. Let $
x_{cm}\equiv\int_{D_1} xd\mu(x)
$ be the mass center of $\mu$ in  $D_1$.  The second moment $Q=Q(\mu)$ of $\mu$ is the symmetric bilinear form  defined by
\begin{align*}
	Q(v,w):=\int_{D_1}[( {x}- {x}_{cm})\cdot v][( {x}- {x}_{cm})\cdot w]d\mu( {x}),\qquad \forall\, v,w\in \R^n.
\end{align*}
Then $Q$ has $n$ eigenvalues (with multiplicity) which can be written as $\lambda_1(\mu)\ge \cdots\ge \lam_n(\mu)$ together with corresponding eigenvectors $v_1(\mu),\cdots,v_n(\mu)$, such that
\begin{align}\label{il}
	Q(v_k)=\int_{D_1}[( {x}- {x}_{cm})\cdot v_k]( {x}- {x}_{cm}) d\mu( {x})=\lam_k v_k.
\end{align}
Describing the eigenvalues by variational method, there holds
\begin{align*}
	\lam_1=\lam_1(\mu):=\max_{|v|^2=1}\int_{D_1}|( {x}- {x}_{cm})\cdot v|^2d\mu( {x}).
\end{align*}
Let $v_1=v_1(\mu)$ be any unit vector achieving such maximum. By induction we have
\begin{align*}
	\lam_{k+1}=\lam_{k+1}(\mu):=\max\left\{\int_{D_1}|( {x}- {x}_{cm})\cdot v|^2d\mu( {x}):{|v|^2=1} ,\ v\cdot v_i=0,\ \text{for all }i\leq k\right\},
\end{align*}
and $v_{k+1}=v_{k+1}(\mu)$ is any unit vector achieving such maximum. Using these eigenvalues and eigenvectors,  it is then not difficult to show that  $V_k= {x}_{cm}+\mbox{span}\{v_1,\cdots,v_k\}$ is a $k$ dimensional affine subspace achieving the minimum in the definition of $\beta_2$;  that is,
\begin{align*}
	\beta_{\mu}^k(0,1)^2=\int_{D_1} d^2( {x},V_k)d\mu( {x})=\lam_{k+1}(\mu)+\cdots+\lam_n(\mu),
\end{align*}
see \cite[Remark 49]{Naber-V-2018}. Using this formula, we can obtain the following estimate, as that of \cite[Proposition 50]{Naber-V-2018}.
\begin{proposition}\label{hj}
	Let  $\mu$ be a probability measure on $D_1$ and $u\in \widehat{H}^{s}(D_{10}, \mathbb{S}^{d-1})$. Let $\lam_k,\ v_k$ be defined  as  above. Then there exists $C(n,s)>0$ such that
	\begin{align*}
		\lam_k\int_{B^+_4\backslash B^+_2}z^a|(v_k,0)\cdot\na u^e(\mathbf{y})|^2d\mathbf{y}\leq C(n,s)\int_{D_1} W_1({x})d\mu({x}), \qquad \text{for all }\,k\ge 1.
	\end{align*}
\end{proposition}

\begin{proof} Without loss of generality we assume  $x_{cm}=0.$ For any $\mathbf{y}=(y,z)\in \R_+^{n+1}$  and $k=1,\cdots,n,$ multiplying $z^{a/2}\na u^e(\mathbf{y})$ on both sides of \eqref{il} yields 	\begin{align}\label{aam} 	
		\lam_kz^{a/2}(v_k,0)\cdot\na u^e(\mathbf{y})=\int_{D_1}z^{a/2}({x}\cdot v_k)(\mathbf{x}\cdot\na u^e(\mathbf{y}))d\mu(x)
			\end{align}
	with $\mathbf{x}=(x,0).$	By the definition of mass center,
	\begin{align*}		\int_{D_1} \mathbf{x}\cdot \mathbf{y}d\mu(x)=\int_{D_1} {x}\cdot {y}d\mu(x)={x}_{cm}\cdot {y}=0.	\end{align*}
	Hence by  \eqref{aam}, 	
	\begin{align*}	
		\lam_kz^{a/2}(v_k,0)\cdot\na u^e(\mathbf{y})=\int_{D_1}z^{a/2}({x}\cdot v_k)((\mathbf{x}-\mathbf{y})\cdot\na u^e(\mathbf{y}))d\mu(x).
			\end{align*}
	By H\"{o}lder's inequality we have
	\begin{align*}		\lam_k^2z^a|(v_k,0)\cdot\na u^e(\mathbf{y})|^2\leq\lam_k\int_{D_1} z^a|((\mathbf{x}-\mathbf{y})\cdot\na u^e(\mathbf{y}))|^2d\mu(x).	\end{align*}
	Without loss of generality we can assume $\lam_k>0$, otherwise there is nothing to prove. Direct computation gives
	\begin{align*}		\lam_k\int_{B^+_4\backslash B^+_2}z^a|\na u^e(\mathbf{y})\cdot (v_k,0)|^2d\mathbf{y}&\leq \int_{D_1}\int_{B^+_4\backslash B^+_2}z^a|(\mathbf{x}-\mathbf{y})\cdot\na u^e(\mathbf{y})|^2d\mathbf{y}d\mu(x)\\	
		&\leq C(n,s) \int_{D_1}\int_{B^+_8(\mathbf{x})\backslash B^+_1(\mathbf{x})}z^a\frac{|(\mathbf{x}-\mathbf{y})\cdot\na u^e(\mathbf{y})|^2}{|\mathbf{x}-\mathbf{y}|^{n+2-2s}} d\mathbf{y}d\mu(x)\\	
		&\leq C(n,s)\int_{D_1} W_1(x)d\mu(x).	\end{align*}
	The proof is complete.\end{proof}

Next we estimate the left hand side of the inequality in the above proposition.
\begin{lemma}\label{lem: non higher order sym}
	For any $\ep>0$, there exists  $\de_8=\de_8(n,\La,s,\ep)>0$ satisfying the following property. For any stationary $s$-harmonic map $u\in \widehat{H}^{s}_\La(D_{10},\mathbb{S}^{d-1})$, if $u$ is  $(0,\de_8)$-symmetric but not   $(k+1,\ep)$-symmetric  on $B^+_1$, then 	
	\begin{align}\label{81}
		\int_{B^+_4\backslash B^+_2}z^a|P\cdot\nabla u^e(\mathbf{y})|^2d\mathbf{y}>\delta_8
	\end{align}
	for every $k+1$ dimensional subspace $P\subset\R^n$.
\end{lemma}
\begin{proof}
	Assume by contradiction that there is a sequence of stationary $s$-harmonic maps $u_i\in \widehat{H}^{s}_\La(D_{10},\mathbb{S}^{d-1})$ such that $u_i$ is  $(0,i^{-1})$-symmetric but not  $(k+1,\ep)$-symmetric  on $B^+_1$  for some $\ep>0$. Moreover, after an orthogonal transformation, there is a  $k+1$ dimensional subspace $P$ such that
	\begin{align}\label{999}
		\int_{B^+_4\backslash B^+_2}z^a|P\cdot\nabla u^e_i(\mathbf{y})|^2d\mathbf{y}\leq i^{-1},\quad i\geq 1.
	\end{align}
	Extracting a subsequence, we can assume that
    \[
    \begin{cases}
    u_i\wto v \quad &\text{in}\,\, \widehat{H}^s(D_{10},\mathbb{R}^{d}),\\
    u_i^e\wto v^e \quad &\text{in}\,\, H^1(B_{10}^+,|z|^ad\mathbf{x}).
    \end{cases}
    \]
    By the compactness theorem \ref{thm: compactness}, we also have strong convergence in $B_4^+$ in the above. It follows that
    \begin{align*}
    \int_{B^+_4\backslash B^+_2}z^a|P\cdot\nabla v^e(\mathbf{y})|^2d\mathbf{y}=0.
    \end{align*}
    Moreover,  using the unique continuation property  and the assumption in above, we know that $v^e$ is $0$-homogeneous in $\mathbb{R}_+^{n+1}$.  	Now we set  $w(\mathbf{y}):=v^e(2\mathbf{y})$, and then
	\[
	\begin{aligned}
	\int_{B^+_2\backslash B^+_1}z^a|P\cdot\nabla v^e(\mathbf{y})|^2d\mathbf{y}&=
	\int_{B^+_2\backslash B^+_1}z^a|P\cdot\nabla w(\mathbf{y})|^2d\mathbf{y}\\
	&=2^{1-n-a}\int_{B^+_4\backslash B^+_2}z^a|P\cdot\nabla v^e(\mathbf{y})|^2d\mathbf{y}=0,
	\end{aligned}
	\]
	where we use the fact that  $w(\mathbf{y})=v^e(\mathbf{y})$ due to the 0-homogeneity of $v^e$. Accordingly, we get
	$\int_{B_4^{+}}z^a|P\cdot\nabla v^e(\mathbf{y})|^2d\mathbf{y}=0.$
	Thus $v$  is  $(k+1)$-symmetric on $B^+_1$. Consequently $u_i$ is $(k+1,\ep)$-symmetric for $i\gg 1$, since $u^e_i$ converges strongly to $v^e$ in $L^{\gamma_0}$. We reach a contradiction. The proof is complete.
\end{proof}
Now we can prove Theorem \ref{ji}.
\begin{proof}[Proof of Theorem \ref{ji}]  By scaling, we may assume $\mu(D_1)=1.$ Since $\lambda_k$ is nonincreasing,
	\begin{align}\label{pl}
		\beta_{\mu}^k(0,1)^2=\lam_{k+1}+\cdots+\lam_n\leq(n-k)\lam_{k+1}.
	\end{align}
	Thus it suffices to estimate $\la_{k+1}$.
	For each $j=1,\cdots,k+1$, by Proposition \ref{hj}, we have
	\begin{align*}
		\sum_{j=1}^{k+1}\lam_j\int_{B^+_4\backslash B^+_2}z^a|(v_j,0)\cdot\na u^e(\textbf{z})|^2d\textbf{z}\leq C(k+1)\int_{D_1} W_1(x)d\mu(x).
	\end{align*}
	Let $V^{k+1}=\mbox{span}(v_1,\cdots,v_{k+1})$. Then
	\begin{align*}
		\lam_{k+1}\int_{B^+_4\backslash B^+_2}z^a|V^{k+1}\cdot\na u^e(\textbf{z})|^2d\textbf{z}&\leq\sum_{j=1}^{k+1}\la_j\int_{B^+_4\backslash B^+_2}z^a|\na u^e(\textbf{z})\cdot (v_j,0)|^2d\textbf{z}\\
		&\leq C\int_{D_1} W_1(x)d\mu(x).
	\end{align*}
	Let $\de_8$ be the number defined as in Lemma \ref{lem: non higher order sym} and set $\delta_7\leq\delta_8$. By assumption $u$ is  $(0,\de_7)$-symmetric but not $(k+1,\ep)$-symmetric on $D_8$.  Then by Lemma \ref{lem: non higher order sym}, there holds
	\begin{align*}
		\int_{B^+_4\backslash B^+_2}z^a|V^{k+1}\cdot\na u^e(\textbf{z})|^2d\textbf{z}\geq \de_7,
	\end{align*}
	and so
	\begin{align*}
		\de_7\lam_{k+1}\leq\lam_{k+1}\int_{B^+_4\backslash B^+_2}z^a|V^{k+1}\cdot\na u^e(\textbf{z})|^2d\textbf{z} \leq C\int_{D_1} W_1(x)d\mu(x).
	\end{align*}
	Since $\de_7\leq\de_8$ depends only on $n,\La,s,\ep$, using \eqref{pl} we conclude that
	\begin{align*}
		\beta_{\mu}^k(0,1)^2\leq C(n,N,\La,s,\ep)\int_{D_1} W_1(x)d\mu(x).
	\end{align*}
	This completes the proof.
\end{proof}

\section{Covering lemma}\label{sec:covering lemma}
Following the refined approach of Naber-Valtorta \cite[Section 6.2]{Naber-V-2018},  we prove the following volume estimate for singular set in this section.

\begin{lemma}[Main covering Lemma]\label{lemma:main covering lemma}
	Fix any $\ep>0$ and $0<r<1$. Let $u\in \widehat{H}_{\La}^{s}(D_{10},\mathbb{S}^{d-1})$ be a stationary $s$-harmonic map.
	 Then there exist $\de=\de(n,d,\La,s,\ep)>0$ and $C(n,\delta, \Lambda)>0$ such that the following holds: for any subset
	${S}\su S^k_{\ep,\de r}(u),$ there exists a finite covering of $S\cap D_1$ satisfying
	\begin{align}\label{mcl}
		S\cap D_1\su \bigcup_{x\in\3}D_{r}(x) \quad\text{and}\quad \sum_{x\in\3}r^k\leq C(n,\delta, \Lambda).
	\end{align}
\end{lemma}

Lemma \ref{lemma:main covering lemma} is a direct consequence of the following  iteration lemma.

\begin{lemma}[Energy drop]\label{ti}
	Let $u\in \widehat{H}_\Lambda^{s}(D_{10},\mathbb{S}^{d-1})$ be a stationary $s$-harmonic map.
	Fix any $\ep>0$ and $0<r<R\leq1$. There exist $\de=\de(n,d,\La,s,\ep)>0$ and $C_2(n)$ such that  for any subset
	${S}\su S^k_{\ep,\de r}(u)$,  there exists a finite covering of $S\cap D_R$ satisfying
	\[
	S\cap D_R\su \bigcup_{x\in\3}D_{r_x}(x) \quad\text{with}\quad  r\le  r_x\le R
	\]
	and
	\[\sum_{x\in\3}r_x^k\leq C_2(n)R^k.\]
	Moreover, for each $x\in\3,$ one of the following conditions  is satisfied:
	\begin{itemize}
		\item[(i)] $r_x= r;$
		\item[(ii)] we have the following uniform energy drop property:
		\begin{equation}\label{110}
			\sup_{y\in D_{2r_x}(x)}\boldsymbol{\Theta}_s(u^e,\mathbf{y},4r_x)\leq E-\de,
		\end{equation}
		where
		$$
		E:=\sup_{x\in D_{2R}}\boldsymbol{\Theta}_s(u^e,\mathbf{x},4R).
		$$
	\end{itemize}
\end{lemma}

With the help of Lemma \ref{ti}, we can prove Lemma \ref{lemma:main covering lemma} directly as follows.

\begin{proof}[Proof of Lemma \ref{lemma:main covering lemma}]
We first apply Lemma \ref{ti} to the initial set $S\cap D_1$ (i.e. taking $R=1$ in Lemma \ref{ti}) so as to obtain a finite covering $\3_1$ with
$$
S\cap D_1\su \bigcup_{x\in\3_1}D_{r_x}(x)\quad\text{and}\quad\sum_{x\in \3_1} r_x^k\leq C_2(n),
$$
where $r\le r_x\le 1$ in the above.  Denote
 \[
\3_{1,+}=\{x|x\in \3_1\,\,\text{with}\,\, r_x> r\}\quad\text{and}\quad\3_{1,r}=\{x|x\in \3_1\,\,\text{with}\,\, r_x= r\}.
\]
For each $x\in \3_{1,+}$,  the assertion (ii) of Lemma \ref{ti}  hold. We can apply Lemma \ref{ti} again to $D_{r_x}(x)$ and obtain a covering  satisfying
\[
S\cap D_{r_x}(x)\su \bigcup_{y\in\3_{2,x}}D_{r_y}(y),\quad r\leq r_y<r_x\quad\text{and}\quad\sum_{y\in \3_{2,x}} r_y^k\leq C_2(n)r_x^k.
\]
Moreover, for each $y\in\3_{2,x}$, if $r_y>r$, then the uniform energy drop holds:
	\[
		\sup_{z\in  D_{2r_y}(y)}\boldsymbol{\Theta}_s(u^e,\mathbf{z},4r_y)\leq E-2\de,
\]

Denote all the ball centers from the above two steps by
$$
\3_2=\bigcup_{x\in \3_{1,+}}\3_{2,x}\cup \3_{1,r}.
$$
Then the construction in the above gives
\[
S\cap D_1\su \bigcup_{x\in\3_2}D_{r_x}(x)\quad\text{and}\quad\sum_{x\in \3_2} r_x^k\leq C_2(n)^2.
\]
Note that the energy $E$ defined as in Lemma \ref{ti} satisfies $E\leq C\Lambda$. Hence the process will stop after at most $i_\delta=\lfloor\de^{-1}C\Lambda\rfloor+1$ iteration. Then we finally obtain a covering $\{D_{r}(x)\}_{x\in \3}$ of $S\cap D_1$ such that
	\[
	\sum_{x\in \3}r^k\leq C_2(n)^{i_\delta}.
	\]
Now it is clear that the estimate \eqref{mcl} is true with $C(n,\delta, \Lambda)=C_2(n)^{\lfloor\de^{-1}C\Lambda\rfloor+1}$.
\end{proof}

Therefore, the main task below is to prove Lemma  \ref{ti}. To this end, We first establish the following  auxiliary Lemma  \ref{oooa}.

\subsection{An auxiliary covering lemma}
\begin{lemma}\label{oooa} For any fixed $\ep>0$ and $0<\rho\leq100^{-1}$, there exist constants $\de=\de(n,d,\La,s,\rho,\ep)>0$ and $C_1(n)$ with the following property.
	For any stationary  $s$-harmonic map $u\in \widehat{H}_{\La}^{s}(D_{10},\mathbb{S}^{d-1})$	 and $0<r<R\leq1$,   let
	$${S}\su S^k_{\ep,\de r}(u)\quad \text{and}\quad E:=\sup\limits_{x\in  D_{2R} }\boldsymbol{\Theta}_s(u^e,\mathbf{x},4R), $$ where $\mathbf{x}=(x,0)$.
	Then,	there exists a finite covering of $S\cap D_R$ such that
	\begin{align*}
		S\cap D_R\su \bigcup_{x\in\3}D_{r_x}(x) \quad\text{with}\quad  R\ge r_x\geq r
	\end{align*}
	and
	\begin{align}\label{pla}
		\sum_{x\in\3}r_x^k\leq C_1(n)R^k.
	\end{align}
	Moreover, for each $x\in\3,$ one of the following conditions  is satisfied:
	\begin{itemize}
	\item[(i)] $r_x=r$;
	\item[(ii)] the set of points
	$$F_x:=\{y\in  D_{2r_x}(x):\boldsymbol{\Theta}_s(u^e,\mathbf{y},\rho r_x/10)>E-\de\}$$
	is contained in $D_{\rho r_x/5}(L_x)\cap D_{2r_x}(x),$ where $L_x$ is a $k-1$ dimensional  affine subspace in $\R^n$.
	\end{itemize}
\end{lemma}
Later we will choose $\rho=\rho(n)$ in the proof of Lemma \ref{ti}, see \eqref{eq:  definition of rho}.
\begin{proof}
	By the scaling invariance property of $\boldsymbol{\Theta}_s$, we may assume $R=1$. For cleanness of the presentation, we also assume in the below that
	\begin{align}\label{332}
		r=\rho^{{l}}, \quad \rho=2^{-a},\quad a,\ {l}\in\N.
	\end{align}

	We divide the proof into two main  steps.
	
	\textbf{Step 1. (Inductive covering)} Let $\eta>0$ be a constant to be determined later. Define the initial set
	\[
	F^0=\{y\in  D_2:\boldsymbol{\Theta}_s(u^e,\mathbf{y},\rho/10)\geq E-\de\}.
	\]
	If there is a $k-1$ dimensional affine subspace $L_0$ such that $F^0\su D_{\rho/5}(L_0),$ then we say that $D_1$ is a good ball. In this case, our claim clearly holds.
	
	Otherwise, $D_1$ is called a bad ball, that is, the set $F^0$ $({\rho}/{10})$-effectively spans a  $k$ dimensional affine subspace $V_0\subset \R^n$ .
	In this case, we choose $\de\leq\de_2$ according to Proposition \ref{pro32} so that
	\begin{equation*}
		S\cap D_1\su D_{\rho/5}(V_0).
	\end{equation*}
	We then choose a finite Vitali covering of $S\cap D_1$ by balls $\{D_\rho(x)\}_{x\in\3^1}$ with $\3^1\su V_0\cap D_1$ such that
	$$
	S\cap D_1\su\bigcup_{x\in \3^1}D_{\rho}(x)\cap D_1.
	$$
	Note that by Proposition \ref{5.5}, for all $x\in\3^1\su V_0\cap D_1$, we have
	$$\boldsymbol{\Theta}_s(u^e,\mathbf{x},\rho/10)\geq E-\eta,$$
	as long as $\de\le \de_3$ which depends also on $\eta$. Under the same smallness assumption ($\delta, \eta\leq \delta_4$), Proposition \ref{Non-symmetry}  implies that for each $x\in \3^1$ we have $x\in S^k_{\ep_{\ast},\rho}(u).$	
	
	Now we divide the above covering according the following rules: for each $x\in\3^1$,  set
	\begin{align*}
		F_x^1=\left\{y\in D_{2\rho}(x):\boldsymbol{\Theta}_s(u^e,\mathbf{y},\rho^2/10)\geq E-\de\right\},
	\end{align*} and  let
	\begin{align*}
		\3^1_g:=\left\{x\in\3^1:\ F_x^1 \su D_{ \rho^2/5}(L_x^1) \
		\mbox{\rm for some}\ k-1\,\mbox{\rm  dimensional affine subspace}\ L_x^1\right\}
	\end{align*}
	and
	\begin{align*}
		\3^1_b:=\{x\in&\3^1:\ F_x^1\,\, (\rho^2/10)\mbox{-effectively spans a}\ k\,\mbox{dimensional affine subspace}\ L_x^{1'}\}.
	\end{align*}
	Then we have $	\3^1=	\3^1_g\cup 	\3^1_b$ and
	$$
	S\cap D_1\subset \bigcup_{x\in \3^1_b}D_{\rho}(x)\cup\bigcup_{x\in \3^1_g}D_{\rho}(x).
	$$
	Hence the remaining work is to deal with bad balls.
	
	Consider the set
	\begin{equation}\label{I_1}
	I_1:=(S\cap D_1)\backslash\bigcup_{x\in \3_g^1}D_\rho(x)=(S\cap D_1)\backslash D_\rho(\3_g^1).
	\end{equation}
	Note that this set is contained in $D_\rho(\3_b^1)$.
	For any bad ball $D_\rho(x)$,
	by Proposition \ref{pro32}, we have
	\begin{equation}\label{SL_x}
	S_{\ep,\delta r}^k\cap D_\rho(x)\su D_{ \rho^2/5}(L_x^{1'})
	\end{equation}
	for some $k$ dimensional affine subspace $L_x^{1'}$. We then consider the set for covering of $I_1$
	\begin{equation}\label{A_1}
	A_1:=\bigcup_{x\in \3_b^1}(D_\rho(x)\cap L_x^{1'})\backslash D_{\rho/2}(\3_g^1).
	\end{equation}
	By \eqref{I_1} and \eqref{SL_x}, we have
	\[
	I_1\subset D_{\rho^{2}/5}(A_1).
	\]
	
	All the balls in $\3^1$ are disjoint under $1/5$ scale, a.e. for each $x_1\neq x_2\in\3^1$, we have $D_{r_{x_1}/5}\cap D_{r_{x_2}/5}=\emptyset$.
	Furthermore, with the same reasoning as the initial covering, there holds
	for all $x\in\3^1$ that $\boldsymbol{\Theta}_s(u^e,\mathbf{x},r_x/10)\geq E-\eta$ and
	for all $s\in[r_{x},1]$, $u$ is not $(k+1,\ep/2)$-symmetric on $D_{s}(x)$, that is, $x\in S^k_{\ep_{\ast},r_{x}}(u)$.
	
	Repeating this procedure, we then build a covering of the form
	\[
	S\cap D_1\su \bigcup_{x\in \3^j}D_{r_{x}}(x)=\bigcup_{x\in \3^j_g}D_{r_{x}}(x)\cup\bigcup_{x\in\3^j_b}D_{r_{x}}(x),
	\]
	with $$\3_b^j=\bigcup\limits_{x\in\3^{j-1}_b}\3^j_{x,b}, \qquad \3^j_g=\left(\bigcup\limits_{x\in\3^{j-1}_b}\3^j_{x,g}\right)\cup\3^{j-1}_g$$
	and $\3^j=\3^j_b\cup\3^j_g,$ for $j\geq2$. Notice that
	\[
	r_{x}=\begin{cases}
		\rho^i, & \text{ if } x\in\3^i_g\ba\3^{i-1}_g,i=2,\cdots,j,\\
		\rho^j, & \text{ if } x\in\3^j_b.
	\end{cases}
	\]
	Moreover, $\3^j_{g}$ is defined to be
	\begin{align*}
		\{x\in\3^j:\ r_{x}\geq\rho^j\ \mbox{\rm and}\ F_x^j\su D_{ \rho r^j_{x}/5}(L_x^j) \
		\mbox{\rm for some}\ k-1\,\mbox{\rm  dimensional affine subspace}\ L_x^j\}
	\end{align*}
	with
	\begin{align*}
		F_x^j:=\{y\in  D_{2r_{x}}(x):\boldsymbol{\Theta}_s(u^e,\mathbf{y},\rho r_{x}/10)\geq E-\de\}
	\end{align*}
	and
	\begin{align*}
		\3^j_{b}:=\{x\in\3^j:& r_{x}=\rho^j \ \mbox{\rm and}
		\ F_x^j\ \rho r_{x}/10\mbox{-}\mbox{\rm effectively spans a}\ k\,\mbox{\rm dimensional affine subspace}\ L_x^{j'}\}.
	\end{align*}
	Furthermore, we have
	\begin{itemize}
		\item[(P1)] For all $x\in\3^j,$  $\boldsymbol{\Theta}_s(u^e,\mathbf{x},r_{x}/10)\geq E-\eta.$
		\item[(P2)] For all $x\in\3^j$ and for all $s\in[r_{x},1]$, $u$ is not $(k+1,\ep_\ast)$-symmetric on $D_{s}(x)$, that is, $x\in S^k_{\ep_\ast,r_{x}}(u)$.
		\item[(P3)] For all $x\neq y\in\3^j$, $D_{r_x/5}(x)\cap D_{r_y/5}(y)=\emptyset$.
	\end{itemize}
	
	Taking $j=l$ and $\rho^{l}=r$ in \eqref{332} and $\3=\3^{l}$. Then the covering part of Lemma \ref{oooa} is complete. We next consider the volume estimates.
	\medskip
	
	\textbf{Step 2. (Volume estimates)}
	In order to prove the volume estimate \eqref{pla}, we define a measure
	\[
	\mu:=\omega_k\sum_{x\in\3}r_x^k\de_x.
	\]
	and the associated  measures
	\[
    \mu_t:=\omega_k\sum_{x\in\3_t}r_x^k\de_x,\]
	where $\3_t=\{x\in\3:r_x\leq t\}$ for all $ r\leq t\le 1$. Keep in mind that $\cal{C}\subset D_1$ and $\mu_t$ is supported on the discrete set $\cal{C}_t$ for all $r\le t\le 1$.
	
	Now we want to prove inductively on
	$r_j=2^jr,j=0,1,\cdots,al-2$ (where $r_{al-2}=1/8$ by \eqref{332}) that the constant $C_R(n)$ appearing in Reifenberg Theorem \ref{rre}, we have
	\begin{equation}\label{1200}
		\mu_{r_j}(D_{r_j}(x))\leq C_R{r_j}^k\quad\text{for each}\quad x\in D_1.
	\end{equation}
	Note that the initial step is true for $r_j=r$, since $D_{r/5}(x)\cap D_{r/5}(y)=\emptyset$ for all  $x\neq y\in\3_r$, it follows easily that $\mu_r(D_r(z))\leq c(n)r^k$ for all $z\in D_1$.
	
	Assume now for all $i\leq j$  we have
	\begin{align}\label{1300}
		\mu_{r_i}(D_{r_i}(x))\leq C_R {r_i}^k,\quad \forall\,x\in D_1.
	\end{align}
	We next show \eqref{1300} is true for ${j+1}$ and hence it holds for all $r_i\leq 1/8$ by induction.\\
	{\bf Rough estimate}\,\, We first show a bad estimate that \eqref{1300} holds for $j+1$ with constant $c(n)C_R$.	Indeed, for fix $x\in D_1$, we can split $\mu_{r_{j+1}}(D_{r_{j+1}}(x))$ as
	\[
	\begin{aligned}
		\mu_{r_{j+1}}(D_{r_{j+1}}(x))&=\mu_{r_j}(D_{r_{j+1}}(x))+\omega_k\sum_{y\in D_{r_{j+1}}(x)\cap \3_{r_{j+1}}\backslash \3_{r_j}} r_{j+1}^k\\
		&=\mu_{r_j}(D_{r_{j+1}}(x))+\omega_kr_{j+1}^k\sharp(D_{r_{j+1}}(x)\cap \3_{r_{j+1}}\backslash \3_{r_j}).
	\end{aligned}
	\]
	Take a covering of $D_{r_{j+1}}(x)$ by balls $D_{r_{j}}(x_i)$ such that $D_{r_{j}}(x_i)$ are disjoint. The number of	$x_i$ has a trivial bound $c(n)$, and by induction we have
	\[
	\mu_{r_j}(D_{r_{j+1}}(x))\leq \sum_i \mu_{r_j}(D_{r_{j}}(x_i))\leq c(n)C_Rr_j^k.
	\]
	For the second part of $\mu_{r_{j+1}}$, by (P3), we get immediately
	\[
	\omega_kr_{j+1}^k\sharp(D_{r_{j+1}}(x)\cap \3_{r_{j+1}}\backslash \3_{r_j})\leq c(n)r_{j+1}^k.
	\]
	{\bf Reifenberg estimates}\,\, We will show inductively that we can apply Theorem \ref{rre} to the measures $\mu_{r_{j+1}}(D_{r_{j+1}}(x))$ for each fixed $x\in D_1$. For convenience, we set
	\[
	\bar{\mu}:=\mu_{r_{j+1}}|_{D_{r_{j+1}}(x)}
	\]
	By consider the Jones' number $\beta_{\bar\mu}^k(y,s)$, notice that the integral with respect to $\bar\mu$ is in fact a sum on $x\in \3_{r_{j+1}}$. We let $y\in \3_{r_{j+1}}$ and then:
	\begin{itemize}
		\item If $s\leq r_y/5$, since $D_s(y)\cap \3_{r_{j+1}}=y$, then $\beta_{\bar\mu}^k(y,s)=0$.
		\item If $s> r_y/5$, because  $\boldsymbol{\Theta}_s(u^e,\mathbf{y},8s)\leq E$, by (P1) we have
		$$
		\boldsymbol{\Theta}_s(u^e,\mathbf{y},8s)-\boldsymbol{\Theta}_s(u^e,\mathbf{y},4s)\leq\eta\quad \text{ for each } s\in[ r_{y}/5,1/8].
		$$
		Given $\delta_7(n,N,\Lambda,s,\rho,\ep_\ast)$ as in Theorem \ref{ji}, it follows from Proposition \ref{pro1} that there exists $\eta_0=\eta_0(n,N,\Lambda,s,\rho,\delta_7)$ such that if $\eta\leq\eta_0$, then $u$ is $(0,\de_7)$-symmetric on $D_{8s}(y)$. By (P2), $u$ is not $(k+1,\ep_\ast)$-symmetric on  $D_{8s}(y)$. Hence we can use Theorem \ref{ji} with $r=s$ then obtain:
		\begin{equation}\label{beta W mu}
			\beta_{\bar\mu}^k(y,s)^2\leq C_1 s^{-k}\int_{D_s(y)}W_s(z)d\bar\mu(z).
		\end{equation}
	\end{itemize}
	By the conditions taking above, we can define the more precisely function:
	\[
	\widehat{W}(x,r):=
	\begin{cases}
		W_r(x)\quad &\text{if}\,\,x\in \3_{r_{j+1}}\,\,\text{and}\,\,r\geq r_x/5,\\
		0\quad&\text{other case},
	\end{cases}
	\]
	and then
	\begin{equation}\label{jones est}
		\beta_{\bar\mu}^k(y,s)^2\leq C_1 s^{-k}\int_{D_s(y)}\widehat{W}(z,s)d\bar\mu(z).
	\end{equation}
	
	We now use Theorem \ref{rre} to the measures $\bar\mu$. For all $\omega\in D_{r_{j+1}}(x)$ and $\tau< r_{j+1}$, we can apply \eqref{jones est}, then have
	\begin{equation}\label{beta mu}
		\begin{aligned}
			& \int_{D_\tau(w)} \int_0^\tau \beta_{\bar\mu}^k(y, s)^2 \frac{d s}{s} d \bar\mu(y) \\
			& \quad \leq C_1\int_{D_\tau(w)} \int_0^\tau  s^{-k} \int_{D_s(y)} \widehat{W}(z,s) d \bar\mu(z) \frac{d s}{s} d \bar\mu(y) \\
			& \quad \leq C_1\int_0^\tau  s^{-k} \int_{D_\tau(w)} \int_{D_s(y)} \widehat{W}(z,s) d \bar\mu(z) d \bar\mu(y) \frac{d s}{s} .
		\end{aligned}
	\end{equation}
	On $D_s(y)$, either
	\(\bar{\mu}=\mu_s|_{D_{r_{j+1}}(x)}\),
	or there exists an \(a\in\operatorname{supp}(\mu)\cap D_s(y)\)
	with \(r_a>s\).
	Since \(y\in\operatorname{supp}(\mu)\) as well, by construction we have
	\(y=a=\operatorname{supp}(\mu)\cap D_s(y)\), and
	\(\widehat{W}(y,s)=0\).
	Thus in either case we have
	\begin{equation}\label{beta mu 2}
		\begin{aligned}
			&\int_{D_\tau(w)} \int_0^\tau \beta_{\bar\mu}^k(y, s)^2 \frac{d s}{s} d \bar\mu(y)\\
			&\quad \leq C_1\int_0^\tau  s^{-k} \int_{D_\tau(w)\cap D_{r_{j+1}}(x)} \int_{D_s(y)\cap D_{r_{j+1}}(x)} \widehat{W}(z,s) d \mu_s(z) d \mu_s(y) \frac{d s}{s} .
		\end{aligned}
	\end{equation}
	Switch the integrals in \eqref{beta mu 2}, we have
	\begin{equation}
		\begin{aligned}
			& \quad\int_{D_\tau(w)} \int_0^\tau \beta_{\bar\mu}^k(y, s)^2 \frac{d s}{s} d \bar\mu(y) \\
			&  \leq C_1\int_0^\tau  s^{-k} \int_{D_{s+\tau}(w)\cap D_{r_{j+1}}(x)} \widehat{W}(z,s)\left(\int_{D_s(z)\cap D_{r_{j+1}}(x)} d \mu_s(y)\right) d \mu_s(z) \frac{d s}{s} \\
			& \leq C_1\int_0^\tau  s^{-k} \int_{D_{2 \tau}(w)\cap D_{r_{j+1}}(x)} \widehat{W}(z,s) \mu_s\left(D_s(z)\right) d \mu_s(z) \frac{d s}{s}.
		\end{aligned}
	\end{equation}
	By induction, and the rough estimates, for all $s\in (0,r_{j+1}]$ and $z\in D_1$ we have
	\begin{equation}
		\mu_s(D_s(z))\leq c(n)C_Rs^k.
	\end{equation}
	Thus we obtain
	\begin{equation}\label{beta W}
		\int_{D_\tau(w)} \int_0^\tau \beta_{\bar\mu}^k(y, s)^2 \frac{d s}{s} d \bar\mu(y)\leq c(n) C_1 C_R  \int_{D_{2 \tau}(w)}\int_0^\tau \widehat{W}(z,s)  \frac{d s}{s} d \bar\mu(z)
	\end{equation}
	Observe by monotonicity formula \eqref{1.8} and (P1), we have
	\[
	\begin{aligned}
		\int_0^\tau\widehat{W}(z,s)\frac{ds}{s}=&\int_{r_z/5}^\tau\widehat{W}(z,s)\frac{ds}{s}\leq\int_{r_z/5}^{1/8}\widehat{W}(z,s)\frac{ds}{s}\\
		\leq& c_1\big(\boldsymbol{\Theta}_s(u^e,\mathbf{z},1))-\boldsymbol{\Theta}_s(u^e,\mathbf{z},r_z/5)\big)\leq c\eta.
	\end{aligned}
	\]
	Plugging this inequation into \eqref{beta W}, and using again the induction hypothesis and the rough estimates, we get
	\begin{equation}\label{beta mu2}
		\int_{D_\tau(w)} \int_0^\tau \beta_{\bar\mu}^k(y, s)^2 \frac{d s}{s} d \bar\mu(y)\leq \delta_5^2\tau^k
	\end{equation}
	where we choose $\eta$ small enough such that
	\begin{equation}\label{eta delta5}
		\eta\leq \eta_1(n,N,\Lambda,s,\ep):=\frac{\delta_5^2}{c(n)C_1C_R^2}.
	\end{equation}
	
	Now, we define a measure $\mu_*(y)=r_{j+1}^{-k}\bar\mu(x+r_{j+1}y)$ and use Theorem \ref{rre} to it. For each $z\in D_1$ and $0<r<1$, by applying the change of variable formula for the integral, we have
	\begin{equation}\label{mu*}
		\int_{D_r(z)}\int_0^r\beta_{\mu_*}^k(y,s)^2\frac{ds}{s}d\mu_*(y)=r_{j+1}^{-k}\int_{D_{rr_{j+1}}(x+r_{j+1}z)}\int_0^r\beta_{\mu_*}^k\left(\frac{y-x}{r_{j+1}},s\right)^2\frac{ds}{s}d\bar\mu.
	\end{equation}
	For the Jones' number
	\begin{equation}\label{betamu*}
		\begin{aligned}
			\beta_{\mu_*}^k\left(\frac{y-x}{r_{j+1}},s\right)^2&=\inf_{V\subset \mathbb{R}^n}\int_{D_s(\frac{y-x}{r_{j+1}})}\frac{d^2(z,V)}{s^2}\frac{d\mu_*(z)}{s^k}\\
			&=\inf_{V\subset \mathbb{R}^n}\int_{D_{sr_{j+1}}(y)}\frac{d^2(z,V)}{(sr_{j+1})^2}\frac{d\bar\mu(z)}{(sr_{j+1})^k}\\
			&=\beta_{\bar\mu}^k(y,sr_{j+1})^2
		\end{aligned}
	\end{equation}
	Plugging \eqref{betamu*} into \eqref{mu*}, use \eqref{beta mu2} and \eqref{eta delta5}, we obtain
	\[
	\begin{aligned}
		\int_{D_r(z)}\int_0^r\beta_{\mu_*}^k(y,s)^2\frac{ds}{s}d\mu_*(y)&=r_{j+1}^{-k}\int_{D_{rr_{j+1}}(x+r_{j+1}z)}\int_0^{rr_{j+1}}\beta_{\bar\mu}^k(y,s)^2\frac{ds}{s}d\bar\mu(y)\\
		&\leq r_{j+1}^{-k}\delta_5^2(rr_{j+1})^k=\delta_5^2 r^k.
	\end{aligned}
	\]
	By apply Theorem \ref{rre}, we obtain $\mu_*(D_1)\leq C_R$. Therefore
	\[
	\mu_{r_{j+1}}(D_{r_{j+1}}(x))\leq C_Rr_{j+1}^k,\quad\text{for each}\,\,x\in D_1.
	\]
	Then we have completed the inductive assumption \eqref{1300}.
	
	   Finally, we choose $\delta=\delta(n,N,\Lambda,s,\rho,\ep)>0$ small enough such that (P1),(P2), \eqref{beta W mu} and \eqref{eta delta5} are satisfied with
	   \[
	   \delta\leq \min\{\delta_2,\delta_3,\delta_4\}\quad\text{and}\quad\eta\leq\min\{\eta_0,\eta_1,\delta_4\}.
	   \]
	   Thus \eqref{1200} holds under the assumption of $\delta$ in above. By (P3) and a simply Vitali covering of $D_1$, we have
	   \[
	   \begin{aligned}
	   \mu(D_1)&=\mu_{1/8}(D_1)+\omega_k\sum_{x\in\3,r_x\in(1/8,1]}r_x^k\delta_x(D_1)\\
	   &\leq \sum_{x\in D_1}\mu_{1/8}(D_{1/8}(x))+\omega_k\sum_{x\in\3,r_x\in(1/8,1]}r_x^k\\
	   &\leq c(n)C_R\left(\frac{1}{8}\right)^k+C(n)=: C_1(n).
	   \end{aligned}
	   \]
	    This finishes the proof of Lemma \ref{oooa}
\end{proof}
We remark that by Theorem 1.1 of Mi\'skiewicz \cite{Miskiewicz-18-Finn}, the last smallness assumption \eqref{eta delta5} on $\eta$ is not much necessary. However, this can not remove the smallness condition on $\eta$, since, for instance, we still need $\eta$ to be very small in the formula below \eqref{beta W mu} in order to apply Proposition \ref{pro1}.

\subsection{Proof of Lemma  \ref{ti}}
By scaling invariance assume without loss of generality that
$$r=\rho^l, \qquad R=1.$$
	In view of Lemma \ref{oooa}, it suffices to prove that Lemma \ref{claim} which can viewed as a weak version of Lemma  \ref{ti}.

\begin{lemma}\label{claim}
 There exists a finite  covering of $S\cap D_1$ by
\[
S\cap D_1\su \bigcup_{x\in\3}D_{r_x}(x)=\bigcup_{x\in\3_{r}}D_{r_x}(x)\cup\bigcup_{x\in\3_{+}}D_{r_x}(x),
\]
satisfying
\[\sum_{x\in\3}r_x^k\leq C^\prime_{2}(n)\]
for some constant $C_2^\prime(n)$ depending only on $n$,	where $\3_{r}$ consists of centers of  $r$-balls in  $\3$, and  $\3_{+}$ consists of centers of  balls $D_{r_x}(x)$ with $r_x> r$ and
\[
\begin{aligned}
\sup_{y\in  D_{2r_x}(x)}\boldsymbol{\Theta}_s(u^e,\mathbf{y},\rho r_x/10)\leq E-\de,
\end{aligned}
\]
where $E,\de$ are the constants defined as in Lemma \ref{oooa}.
\end{lemma}
\begin{proof}
 We still apply an induction argument. By superscripts
	$f,b$ we mean \emph{final} and \emph{bad} balls respectively.
	
	\textbf{Step 1. (Recovering of bad balls: the first step)} The  auxiliary Lemma \ref{oooa}  gives a covering $\{D_{r_x}(x)\}_{x\in \3}$ of $S\cap D_{1}$ satisfying
	\[
	S\cap D_{1}\subset\bigcup_{x\in{\cal C}_{r}^{0}}D_{r}(x)\cup\bigcup_{x\in{\cal C}_{+}^{0}}D_{r_{x}}(x),
	\]
	where
	\[
	{\cal C}_{r}^{0}=\{x\in{\cal C}:r_{x}=r\}\quad\text{and}\quad{\cal C}_{+}^{0}=\{x\in{\cal C}:r_{x}>r\},
	\]
	and
	\begin{equation}\label{eq: 6.67}
		\sum_{x\in{\cal C}_{r}^{0}\cup{\cal C}_{+}^{0}}r_{x}^{k}\le C_1(n).
	\end{equation}
	Moreover, for each $x\in{\cal C}_{+}^{0}$, the set
	\[
	F_{x}=\left\{ y\in  D_{2r_{x}}(x):\boldsymbol{\Theta}_s(u^e,\mathbf{y},\rho r_x/10)>E-\de\right\}
	\]
	is contained in a small neighborhood of a ($k-1$)-dimensional affine subspace. To proceed, we only need to refine the part ${\cal C}_{+}^{0}$.

	Let $x\in{\cal C}_{+}^{0}$ so that $r_x>r$. We have two cases:
	
	\textbf{Case 1:}  $\rho r_{x}=r$. In view of the desired covering in the statement of Lemma \ref{ti}, we directly cover $S\cap D_{r_{x}}(x)$
	by a family of balls
	$$\Big\{D_{r_{y}}(y): y\in{\cal C}_{x}^{(1,r)}\subset S\cap D_{r_{x}}(x)\Big\}\quad \text{with}\quad r_{y}=\rho r_{x}=r,$$
	and $\{D_{r_{y}/5}(y)\}_{y\in{\cal C}_{x}^{(1,r)}}$ being pairwise disjoint.  A simple volume comparison argument shows that the cardinality  $\sharp\left({\cal C}_{x}^{(1,r)}\right)\le C(n)\rho^{-n}$. Hence
	\[
	\sum_{y\in{\cal C}_{x}^{(1,r)}}r_{y}^{k}=\sharp\left({\cal C}_{x}^{(1,r)}\right)(\rho r_{x})^{k}\le C(n)\rho^{k-n}r_{x}^{k}=:C_{r}(n,\rho)r_{x}^{k}.
	\]
	Collect all such points $x$ and set
	\[
	{\cal C}^{(1,r)}={\cal C}_{r}^{0}\cup\bigcup_{x\in{\cal C}_{+}^{0},\rho r_{x}=r}{\cal C}_{x}^{(1,r)}.
	\]
	It follows from \eqref{eq: 6.67} that
	\[
	\sum_{y\in{\cal C}^{(1,r)}}r_{y}^{k}=\left(\sum_{x\in{\cal C}_{r}^{0}}+\sum_{x\in{\cal C}_{+}^{0}}\sum_{y\in{\cal C}_{x}^{(1,r)}}\right)r^{k}\le\sum_{x\in{\cal C}_{r}^{0}}r^{k}+C_{r}(n,\rho)\sum_{x\in{\cal C}_{+}^{0}}r_{x}^{k}\le C_{1}(n)C_{r}(n,\rho).
	\]
	
	\textbf{Case 2}:  $\rho r_{x}>r$.
	
	In this case, we  cover $S\cap D_{r_{x}}(x)$ by making  use of the key fact that $F_{x}\subset D_{\rho r_{x}/5}(L)\cap D_{r_{x}}(x)$ for some $k-1$ dimensional affine subspace $L$.
	First we choose a Vitali  covering for the part away from $D_{2\rho r_{x}}(F_{x})$:
	\[
	S\cap D_{r_{x}}(x)\Big\backslash D_{2\rho r_{x}}(F_{x})\subset\bigcup_{y\in{\cal C}_{x}^{(1,f)}}D_{r_{y}}(y)\quad\text{with }r_{y}=\rho r_{x},
	\]
	where ${\cal C}_{x}^{(1,f)}\subset S\cap D_{r_{x}}(x)$. Since $y$ is away from  $F_{x}$,  the energy drop property holds:
	\begin{equation}\label{eq: energy-drop-2}
		\boldsymbol{\Theta}_s\left(u^e, \mathbf{y},r_y/10\right)\le E-\de.
	\end{equation}
	In fact, we also have the almost desired energy drop property:
	\begin{equation}\label{eq: energy-drop-3}
		\sup_{z\in  D_{2r_{y}}(y)} \boldsymbol{\Theta}_s(u^e, \mathbf{z},\rho r_y/10)\le E-\de, \quad \forall\,y\in{\cal C}_{x}^{(1,f)}.
	\end{equation}
	Indeed, for each $z\in S\cap D_{2r_{y}}(y)$,
	the fact $|z-y|\le2r_{y}=2\rho r_{x}$ implies that $z\not\in F_{x}$.
	Note also that $ D_{2r_{y}}(y)\subset  D_{2r_{x}}(x)$ since
	$\rho=\rho(n)$ will be chosen less than $1/100$. Hence $z\in  D_{2r_{x}}(x)$.
	Consequently, the definition of $F_x$ implies that
	\[
	\boldsymbol{\Theta}_s(u^e, \mathbf{z},\rho r_y/10)\le\boldsymbol{\Theta}_s(u^e, \mathbf{z},\rho r_x/10)\le E-\de.
	\]
	This proves \eqref{eq: energy-drop-3}. Due to this energy drop  property,	these balls $\{D_{r_{y}}(y): y\in {\cal C}_{x}^{(1,f)}\}$ will be part of the ``final'' balls. Moreover, the number $\sharp\left({\cal C}_{x}^{(1,f)}\right)$ is bounded from above by a constant $C(n)\rho^{-n}$. This implies that
	\[
	\sum_{y\in{\cal C}_{x}^{(1,f)}}r_{y}^{k}=\sharp\left( {\cal C}_{x}^{(1,f)}\right)(\rho r_{x})^{k}\le C(n)\rho^{k-n}r_{x}^{k}=C_{f}(n,\rho)r_{x}^{k}.
	\]
	
	Collecting all the subset ${\cal C}_{x}^{(1,f)}$ from the above,  we
	obtain the first generation of final balls:
	\[
	{\cal C}^{(1,f)}=\bigcup_{x\in{\cal C}_{+}^{0},\rho r_{x}>r}{\cal C}_{x}^{(1,f)}.
	\]
	together with the volume estimate  from \eqref{eq: 6.67}  that
	\[
	\sum_{y\in{\cal C}^{(1,f)}}r_{y}^{k}=\sum_{x\in{\cal C}_{+}^{0}}\sum_{y\in{\cal C}_{x}^{(1,f)}}(r_y)^{k}\le C_{f}(n,\rho)\sum_{x\in{\cal C}_{+}^{0}}r_{x}^{k}\le C_{f}(n,\rho)C_{1}(n).
	\]
	
	For the remaining part of $S\cap D_{r_{x}}(x)$, we first choose a Vitali covering:
	\[
	S\cap D_{r_{x}}(x)\cap D_{2\rho r_{x}}(F_{x})\subset\bigcup_{y\in{\cal C}_{x}^{1,b}}D_{r_{y}}(y)\quad\text{with }r_{y}=\rho r_{x},
	\]
	where ``$b$'' means bad balls, on which $r_y>r$ and
	the	energy drop property \eqref{eq: energy-drop-2}   can not be  determined as well.
	However, since $F_{x}\subset D_{\rho r_{x}/5}(L)\cap D_{2r_{x}}(x)$
	for some ($k-1$)-dimensional affine subspace $L$, we have a better cardinality estimate by volume comparison 	argument:
	\[
	\sharp\left(\3_{x}^{(1,b)}\right)\le C_b(n)\rho^{1-k}.
	\]
	This means that relatively there are  quite few  bad balls are. As a result,
	\[
	\sum_{y\in{\cal C}_{x}^{(1,b)}}(r_{y})^{k}=\sharp\left(\3_{x}^{(1,b)}\right)(\rho r_{x})^{k}  \le C_{b}(n)\cdot \rho \cdot r_{x}^{k}.
	\]
	Collect all the bad balls together to get the first generation of bad balls
	\[
	{\cal C}^{(1,b)}=\bigcup_{x\in{\cal C}_{+}^{0},\rho r_{x}>r}{\cal C}_{x}^{(1,b)},
	\]
	together with the volume estimate by \eqref{eq: 6.67}
	\[
	\sum_{y\in{\cal C}^{(1,b)}}r_{y}^{k}=\sum_{x\in{\cal C}_{+}^{0}}\sum_{y\in{\cal C}_{x}^{(1,b)}}(r_{y})^{k}\le C_{b}(n)\rho\sum_{x\in{\cal C}_{+}^{0}}r_{x}^{k}\le C_{1}(n)C_{b}(n)\rho.
	\]
	Now we choose $\rho=\rho(n)$ by requiring
	\begin{equation}\label{eq:  definition of rho}
		0<\rho<\min\left\{ 100^{-1},\frac{1}{2C_{1}(n)C_{b}(n)}\right\}
	\end{equation}
	such that
	\[
	\sum_{y\in{\cal C}^{(1,b)}}(\rho r_{x})^{k}\le 1/2.
	\]
	Note that $r<r_{y}=\rho r_x\le\rho$ for all $y\in{\cal C}^{(1,b)}$.

	\textbf{Step 2. (Induction step)}  To proceed, denote
	\[
	C^\prime_{2}(n)=2C_{1}(n)\left(C_{r}(n,\rho)+C_{f}(n,\rho)\right),
	\]
	such that
	\[
	\sum_{y\in{\cal C}^{(1,r)}\cup{\cal C}^{(1,f)}}r_{y}^{k} \le(C_{f}(n,\rho)+C_{r}(n,\rho))\sum_{x\in{\cal C}}r_{x}^{k}\le\frac{1}{2}C^\prime_{2}(n)
	\] and that
	\[
	\sum_{y\in{\cal C}^{(1,b)}}(\rho r_{x})^{k}\le \frac12\le   \frac12C^\prime_{2}(n).
	\]	
	Our aim is to derive the following
	covering: for each $1\le i\le l$, there holds
	\[
	S\cap D_{1}\subset\bigcup_{x\in{\cal C}^{(i,r)}}D_{r_x}(x)\cup\bigcup_{x\in{\cal C}^{(i,f)}}D_{r_{x}}(x)\cup\bigcup_{x\in{\cal C}^{(i,b)}}D_{r_{x}}(x).
	\]
	Moreover,  the following properties hold:
	\begin{itemize}
	\item[(1)] For each $x\in{\cal C}^{(i,r)}$, $r_{x}=r$;
	\item[(2)] For each $x\in{\cal C}^{(i,f)}$, the energy drop condition holds:
	\begin{equation}\label{eq: energy-drop-4}
		\sup_{z\in  D_{2r_{x}}(x)} \boldsymbol{\Theta}_s(u^e, \mathbf{z},\rho r_x/10)\le E-\de;
	\end{equation}
	these balls will be called final balls;
	\item[(3)] For each $x\in{\cal C}^{(i,b)}$, we have $r<r_{x}\le\rho^{i}$.
	On these ``bad'' balls, none of the above two  conditions is
	verified.
	\item[(4)] There holds
	\begin{equation}\label{eq: 6.70}
		\sum_{y\in{\cal C}^{(i,r)}\cup{\cal C}^{(i,f)}}r_{y}^{k}\le\left(\sum_{j=1}^{i}2^{-j}\right)C^\prime_{2}(n)\quad \text{and}\quad \sum_{y\in{\cal C}^{(i,b)}}r_{y}^{k}\le2^{-i}\le2^{-i}C^\prime_{2}(n).
	\end{equation}
    \end{itemize}
	We have proved the above induction for $i=1$ in Step 1. 	Suppose now
	it holds for some $1\le i<l$. We need to prove that it holds for
	$i+1$.		It is clear that we only need to decompose those bad balls.
	
	Fix a point $x\in{\cal C}^{(1,b)}$, i.e., $D_{r_{x}}(x)$ is a bad
	ball and $r<r_{x}\le\rho^{i}$.
	
	\textbf{Case 1: $\rho r_{x}=r$.} Then we cover it simply by balls
	of radius $\rho r_{x}$ as in step 1, and obtain a covering of $r$-balls
	whose centers lie in ${\cal C}_{x}^{(i+1,r)}$, together with the estimate
	\[
	\sum_{y\in{\cal C}_{x}^{(i+1,r)}.}(\rho r_{x})^{k}\le C(n)\rho^{k-n}r_{x}^{k}=C_{r}(n,\rho)r_{x}^{k}.
	\]
	
	\textbf{Case 2: $\rho r_{x}>r$.} In this case we apply the previous lemma
	to get a covering $\{D_{r_{y}}(y)\}_{y\in{\cal C}^{x}}$ such that
	\[
	S\cap D_{r_{x}}(x)\subset\bigcup_{y\in{\cal C}_{r}^{x}}D_{r}(y)\cup\bigcup_{y\in{\cal C}_{+}^{x}}D_{r_{y}}(y)\qquad\text{with }r_{y}\ge r
	\]
	and
	\[
	\sum_{y\in{\cal C}_{r}^{x}\cup{\cal C}_{+}^{x}} r_{y}^{k}\le C_{1}(n)r_{x}^{k}.
	\]
	Moreover, for each $y\in{\cal C}_{+}^{x}$, there is a $k-1$ dimensional
	affine subspace $L_{y}$ such that
	\[
	F_{y}\equiv\left\{ z\in  D_{r_{y}}(y):\boldsymbol{\Theta}_s(u^e,\mathbf{z},\rho r_y/10)\ge E-\de\right\} \subset D_{\rho r_{y}/5}(L_{y})\cap D_{r_{y}}(y).
	\]
	Of course we will reserve ${\cal C}_{r}^{x}$ as part of ${\cal C}^{(i+1,r)}$.
	Thus below we assume that $y\in{\cal C}_{+}^{x}$. The method is totally similar.
	
	\textbf{Case 2.1: $\rho r_{y}=r$.} Then as that of Case 1, we get
	a simple covering of at most $C(n)\rho^{-n}$ $r$-balls $\{D_{r_z}(z)\}_{z\in{\cal C}_{y}^{(i+1,r)}}$
	of $D_{r_{y}}(y)$ with $r_z=r$. So we define
	\[
	{\cal C}^{(i+1,r)}={\cal C}^{(i,r)}\cup\bigcup_{x\in{\cal C}^{(i,b)},\rho r_{x}=r}{\cal C}_{x}^{(i+1,r)}\cup\bigcup_{x\in{\cal C}^{(i,b)},\rho r_{x}>r}\left({\cal C}_{r}^{x}\cup\bigcup_{y\in{\cal C}_{+}^{x},\rho r_{y}=r}{\cal C}_{y}^{(i+1,r)}\right).
	\]
	This shows that how much more is ${\cal C}^{(i+1,r)}$ than that of
	${\cal C}^{(i,r)}$. The newly generated $r$-balls have measure
	\[
	\begin{aligned} & \sum_{x\in{\cal C}^{(i,b)},\rho r_{x}=r}\left(\sum_{z\in{\cal C}_{x}^{(i+1,r)}}r_{z}^{k}\right)
		\qquad+\sum_{x\in{\cal C}^{(i,b)},\rho r_{x}>r}\left(\sum_{y\in{\cal C}_{r}^{x}}r_{y}^{k}+\sum_{y\in{\cal C}_{+}^{x}}\sum_{z\in{\cal C}_{y}^{(i+1,r)}}r_{z}^{k}\right)\\
		& \qquad \le C_{r}(n,\rho)\left(\sum_{x\in{\cal C}^{(i,b)},\rho r_{x}=r}r_{x}^{k}+\sum_{x\in{\cal C}^{(i,b)},\rho r_{x}>r}\left(\sum_{y\in{\cal C}^{x}}r_{y}^{k}\right)\right)\\
		&  \qquad \le C_{r}(n,\rho)C_{1}(n)\sum_{x\in{\cal C}^{(i,b)}}r_{x}^{k}\\
		&  \qquad \le2^{-i}C_{r}(n,\rho)C_{1}(n).
	\end{aligned}
	\]
	Thus the total measure of $r$-balls in the $(i+1)$-th generation is estimated by
	\[
	\begin{aligned}\sum_{z\in{\cal C}^{(i+1,r)}}r_{z}^{k}
		& \le\sum_{z\in{\cal C}^{(i,r)}}r_{z}^{k}+2^{-i}C_{r}(n,\rho)C_{1}(n).
	\end{aligned}
	\]
	
	\textbf{Case 2.2: $\rho r_{y}>r$.} We cover $D_{r_{y}}(y)$ by
	\[
	\begin{aligned} & S\cap D_{r_{y}}(y)\Big\backslash D_{2\rho r_{y}}(F_{y})\subset\bigcup_{z\in{\cal C}_{y}^{(i+1,f)}}D_{r_{z}}(z),\\
		& S\cap D_{r_{y}}(y)\cap D_{2\rho r_{y}}(F_{y})\subset\bigcup_{z\in{\cal C}_{y}^{(i+1,b)}}D_{r_{z}}(z),
	\end{aligned}
	\quad\text{with }r_{z}=\rho r_{y}.
	\]
	As that of \eqref{eq: energy-drop-3},  each ball $D_{r_{z}}(z)$ in the first cover ${\cal C}_{y}^{(i+1,f)}\subset S\cap D_{r_{y}}(y)$ satisfies
	the energy drop property \eqref{eq: energy-drop-4} in the inductive assumption. Moreover,
	\[
	\sum_{z\in{\cal C}_{y}^{(i+1,f)}}r_{z}^{k}\le C_{f}(n,\rho)r_{y}^{k}\quad\text{and}\quad \sum_{z\in{\cal C}_{y}^{(i+1,b)}}r_{z}^{k}\le C_{b}(n)\cdot \rho\cdot r_{y}^{k}.
	\]
	Note that in this case no $r$-balls occur anymore. Set
	\[
	{\cal C}^{(i+1,f)}={\cal C}^{(i,f)}\cup\bigcup_{x\in{\cal C}^{(i,b)}}\bigcup_{y\in{\cal C}_{+}^{x}}{\cal C}_{y}^{(i+1,f)},\quad{\cal C}^{(i+1,b)}=\bigcup_{x\in{\cal C}^{(i,b)}}\bigcup_{y\in{\cal C}_{+}^{x}}{\cal C}_{y}^{(i+1,b)}.
	\]
	We have
	\[
	\begin{aligned}\sum_{z\in{\cal C}^{(i+1,f)}}r_{z}^{k} & =\sum_{z\in{\cal C}^{(i,f)}}r_{z}^{k}+\sum_{x\in{\cal C}^{(i,b)}}\sum_{y\in{\cal C}_{+}^{x}}\sum_{z\in{\cal C}_{y}^{(i+1,f)}}r_{z}^{k}\\
		& \le\sum_{z\in{\cal C}^{(i,f)}}r_{z}^{k}+C_{f}(n,\rho)\sum_{x\in{\cal C}^{(i,b)}}\sum_{y\in{\cal C}_{+}^{x}}r_{y}^{k}\\
		& \le\sum_{z\in{\cal C}^{(i,f)}}r_{z}^{k}+C_{f}(n,\rho)C_{1}(n)\sum_{x\in{\cal C}^{(i,b)}}r_{x}^{k}\\
		& \le\sum_{z\in{\cal C}^{(i,f)}}r_{z}^{k}+2^{-i}C_{f}(n,\rho)C_{1}(n).
	\end{aligned}
	\]
	Recall that $C^\prime_{2}(n)=2\left(C_{r}(n,\rho)+C_{f}(n,\rho)\right)C_{1}(n)$.
	Hence by the inductive assumption we obtain
	\[
	\left(\sum_{z\in{\cal C}^{(i+1,r)}}+\sum_{z\in{\cal C}^{(i+1,f)}}\right)r_{z}^{k}\le\left(\sum_{z\in{\cal C}^{(i,r)}}+\sum_{z\in{\cal C}^{(i,f)}}\right)r_{z}^{k}+2^{-i-1}C^\prime_{2}(n)\le C^\prime_{2}(n)\sum_{j=1}^{i+1}2^{-j},
	\]
	and also by the choice of $\rho$, we obtain
	\[
	\begin{aligned}\sum_{z\in{\cal C}^{(i+1,b)}}r_{z}^{k}= & \sum_{x\in{\cal C}^{(i,b)}}\sum_{y\in{\cal C}_{+}^{x}}\sum_{z\in{\cal C}_{y}^{(i+1,b)}}r_{z}^{k}\le\sum_{x\in{\cal C}^{(i,b)}}\sum_{y\in{\cal C}_{+}^{x}}C_{b}(n)\rho r_{y}^{k}\\
		& \le\sum_{x\in{\cal C}^{(i,b)}}C_{b}(n)\rho C_{1}(n)r_{y}^{k}\le2^{-1}\sum_{x\in{\cal C}^{(i,b)}}r_{y}^{k}\le2^{-i-1}.
	\end{aligned}
	\]
	This proves \eqref{eq: 6.70} for $i+1$.
	
	Finally,  note that for every $z\in{\cal C}^{(i+1,b)}$, there exists $x\in{\cal C}^{(i,b)}$ such that
	$$r<r_{z}=\rho r_{y}\le\rho r_{x}\le\rho^{i+1}, $$
	since $r_{x}\le\rho^{i}$ by inductive	 assumption.   The proof of the induction is complete.
	
	Since $r=\rho^{l}$, the above procedure will stop at $i=l$.	The proof is complete.
\end{proof}

Lemma  \ref{ti} will follow by  decomposing the  balls in $\3_{+}$ into smaller ones in a controllable way.

\begin{proof}[Proof of Lemma  \ref{ti}]
So suppose the above claim holds. Fix a point $x\in\3_{+}$.  If $\rho^2 r_x\le r$, then directly choose a Vitali covering of at most $C(n,\rho)$ balls of $\{D_r(y)\}_{y\in \mathcal{C}^x_r}$ to cover $S\cap D_{r_x}(x)$. The volume   is well controlled:
\[
\sum_{y\in \mathcal{C}^x_r} r_y^k\le C(n,\rho) r_x^k.
\]
In the case $\rho^2 r_x> r$, we cover $S\cap D_{r_x}(x)$ by a Vitali covering of at most $C(n,\rho)$ balls of $\{D_{r_y}(y)\}_{y\in \mathcal{C}^x_+}$ with $r_y=\rho^2 r_x$. Then, we have the following inclusion relationship
\[
S\cap D_{r_x}(x)\subset   \bigcup_{y\in \mathcal{C}^x_+}D_{r_y}(y)\quad\text{and}\quad D_{2r_y}(y)\subset D_{2r_x}(x).
\]
This makes sure that for each $z\in D_{2r_y}(y)$, there holds
\[D_{4r_y}(z)=D_{4\rho^2 r_x}(z)\subset D_{\rho r_x/10}(z). \]
Consequently, by the monotonicity of $\boldsymbol{\Theta}_s$, we deduce, for each $y\in \mathcal{C}^x_+$, that
\[
\begin{aligned}
\sup_{z\in  D_{2r_y}(y)}\boldsymbol{\Theta}_s(u^e,\mathbf{z},4r_y)\leq E-\de.
\end{aligned}
\]
With the help of the above decomposition, we obtain the covering
\[S\cap D_1\subset\bigcup_{x\in\3_{r}}D_{r}(x)\cup\bigcup_{x\in\3_{+}, \rho^2 r_x\le r}\bigcup_{ y\in \mathcal{C}^x_r}D_{r}(y)\cup\bigcup_{x\in\3_{+}, \rho^2 r_x> r}\bigcup_{y\in  \mathcal{C}^x_+}D_{r_y}(y),
\]
together with the volume estimate
\[
\sum_{x\in \mathcal{C}_r}r_x^k + \sum_{x\in \3_+}\sum_{y\in \mathcal{C}^x_r\cup \mathcal{C}^x_+ }r_y^k\leq C(n,\rho)\sum_{x\in \3}r_x^k\leq C(n,\rho)C_2^\prime(n).
\]
This completes the proof of Lemma \ref{ti} upon taking 	$C_2(n)=C(n,\rho)C_2^\prime(n)$.
\end{proof}

Therefore it remains to prove the above claim.

\section{Proof of Theorem \ref{thm80}: Volume estimates and rectifiability}
This theorem follows from the main covering lemma  \ref{lemma:main covering lemma} and the rectifiable-Reifenberg theorem mentioned in the previous sections.
\begin{proof}[Proof of Theorem \ref{thm80}] Let $\de=\delta(n,d,\Lambda,s,\ep)>0$ be defined as in Lemma \ref{lemma:main covering lemma}. It follows that, for all $0<r<\delta$, let $r_0=r/\delta$, we get
	\begin{align}\label{z3}
	 \mbox{Vol}\Big(D_r\Big(S^k_{\ep,r}(u)\cap D_1\Big)\Big)=\mbox{Vol}\Big(D_{\de r_0}\Big(S^k_{\ep,\delta r_0}(u)\cap D_1\Big)\Big) \leq \sum_{x\in\3} (2r_0)^{n}\leq C_\ep r^{n-k},
	\end{align}
	where $C_{\ep}=2^n \delta^{k-n}C(n,\delta, \Lambda)$. This proves the desired volume estimate \eqref{J} for all $0<r<\de$. In the other hand, for any $\de\le r<1$ we have also
	\[\label{z3}
	\mbox{Vol}\Big(D_r\Big(S^k_{\ep,r}(u)\cap D_1\Big)\Big) \le \mbox{Vol}(D_2) \le C(n) \left(\frac{r}{\de} \right)^{n-k}\leq C_\ep r^{n-k}.\]
	The proof of  \eqref{J} is complete. The volume estimate \eqref{Jk} follows from  \eqref{J} by noting that $S^k_{\ep}(u)\subset S^k_{\ep, r}(u)$ for any $r>0$.

	In order to prove  the rectifiability of $S^k(u)$, it is sufficient to prove the rectifiability of $S^k_\ep(u)$ for each $\ep>0$, as $S^k(u)=\bigcup_{i\ge 1}S^k_{1/i}(u)$.
	
	By the volume estimate \eqref{Jk}, we have $\HH^k(S^k_\ep(u)\cap D_1)\leq C_\ep.$ Applying the same estimates on $D_r(x)$ with $x\in D_1$ and $r\leq 1$ gives the Alhfors upper bound estimate
	\begin{equation}\label{jkj}
		\HH^k(S^k_\ep(u)\cap D_r(x))\leq C_\ep r^k.
	\end{equation}
	
	Let $S\su S_\ep^k(u)\cap D_1$ be an arbitrary measurable subset with $\HH^k(S)>0$. Set
	$$g(x,r)=\boldsymbol{\Theta}_s(u^e,\mathbf{x},r)-\Xi(u,x),\quad\forall\, x\in D_1\text{ and }\ r\leq 1,$$
	where $\Xi(u,x)=\lim_{r\to 0}\boldsymbol{\Theta}_s(u^e,\mathbf{x},r)$. Since $u\in \widehat{H}_{\La}^{s}(\Omega,\mathbb{S}^{d-1})$, we know that $g$ is uniformly bounded and converges to zero everywhere,  the dominated convergence theorem implies that for each $\de>0,$ there exists $\overline{r}>0$ such that
	\begin{equation*}\label{139}
		\medint_S g(x,10\overline{r})d\HH^k(x)\leq\de^2.
	\end{equation*}
	So we can find a measurable subset $E\su S$ with $\HH^k(E)\leq\de\HH^k(S)$ and $g(x,10\overline{r})\leq\de$ for all $x\in F:=S\ba E$. Cover $F$ by a finite number of balls $D_{\overline{r}}(x_i)$ centered on $F.$ Rescaling if necessary, we may assume that $D_{\overline{r}}(x_i)=D_1$. Then $g(x,10)\leq\de$ for $x\in F.$ Similar to \eqref{beta W mu}, choose $\de$ sufficiently small such that $u$ is $(0,\de_7)$ symmetric in $D_{10}$. Theorem \ref{ji} implies
	\begin{equation*}
		\beta_{\HH^k|_F}^k(z,s)^2\leq C_1s^{-k}\int_{D_s(z)}W_s(y)d\HH^k|_F(y) \quad \text{for all }z\in F,\ s\leq1.
	\end{equation*}
	Integrating the previous estimate with respect to $z$ and using \eqref{jkj} yield that for all  $x\in D_1$ and $s\leq r\leq1,$ there holds
	\begin{align*}
		\int_{D_r(x)}\beta_{\HH^k|_F}^k(z,s)^2d\HH^k|_F(z)&\leq C_1s^{-k}\int_{D_r(x)}\int_{D_s(z)}W_s(y)d\HH^k|_F(y)d\HH^k|_F(z)\non\\
		&\leq C_1s^{-k}\int_{D_{r+s}(x)}W_s(y)\int_{D_s(y)}d\HH^k|_F(z) d\HH^k|_F(y)\\
		&\leq
		C_1C_\ep\int_{D_{r+s}(x)}W_s(y) d\HH^k|_F(y).
	\end{align*}
	Integrating again with respect to $s$, similar to \eqref{beta mu2}, we obtain that for all $x\in D_1$ and $r\leq1,$
	\begin{align*}
		\int_{D_r(x)}\int_0^r\beta_{\HH^k|_F}^k(z,s)^2\frac{ds}{s}d\HH^k|_F(z)&\leq C_1C_\ep\int_{D_{2r}(x)}\left(\boldsymbol{\Theta}_s(u^e,\mathbf{y},10r)-\Xi(u,y)\right)d\HH^k|_F(y)\non\\
		&\leq c(n)C_1C_\ep^2\de r^k.
	\end{align*}
	Choosing
	$$\de\le\frac{\de_6^2}{c(n)C_1C_\ep^2},$$
	we deduce from Theorem \ref{rree} that $F\cap D_1$ is $k$-rectifiable.
	
	Repeating the above argument with \(E\) in place of the initial set \(S\), we could find another measurable set $E_1\subset E$ with $\HH^k(E_1)\leq \delta \HH^k(E)$, and that $F_1:=E\ba E_1$ is $k$-rectifiable. Continuing this process, we eventually conclude that $S$ is $k$-rectifiable.
	
	The proofs of final assertions are similar to that of \cite{Naber-V-2017} and are thus omitted.
\end{proof}

\section{Proof of Theorem \ref{thm84}}\label{sec: proof of main results}

The proof of Theorem \ref{thm84} follows the same idea of Cheeger and Naber \cite{Cheeger-Naber-2013-CPAM}. It also relies on the following $\ep$-regularity theorem.

\begin{theorem}(\cite[Theorem 4.1]{He-X-Z-24})\label{thm: sym-to-reg} Let  $u\in \widehat{H}_\Lambda^s(D_{2},\mathbb{S}^{d-1})$ be a stationary $s$-harmonic map if $1/2<s<1$, or a minimizing $s$-harmonic map if $0<s\le 1/2$. There exists
	a constant $\de(n,\La,s)>0$ such that, if $u$ is $(n-1,\de)$-symmetric on $D_{1}$, then
	\[
	r_{u}(0)\ge\ka_{2},
	\]
	where $\ka_{2}=\ka_{2}(n,s)>0$	is the constant defined as in Theorem \ref{thm: partial regularity of MPS}.
\end{theorem}

Now we are ready to prove Theorem \ref{thm84}.

\begin{proof}[Proof of Theorem \ref{thm84}] Let $r<\kappa_2$ and $\mathcal{B}_r(u)=\{x\in D_1: r_u(x)< r\}$. A simple scaling argument together with the application of Theorem  \ref{thm: sym-to-reg} implies that
	\[
	\mathcal{B}_r(u)\subset \mathcal{S}_{\delta, r/\kappa_2}^{n-2}.
	\]
	Then Theorem \ref{thm80}  yields
	\[
	\operatorname{Vol}\left(D_r\left(\mathcal{B}_r(u)\right) \cap D_1\right) \leq \operatorname{Vol}\left(D_r\left(\mathcal{S}_{\delta,  r/\kappa_2}^{n-2}\right) \cap D_1\right) \leq C(s, n, \Lambda,d) r^2.
	\]
	
		To estimate $\na u$, simply observe that by Definition \ref{pppv} there holds
	\begin{align*}
	\{x\in D_1(0):|\nabla u|>r^{-1}\}\su\{x\in D_1(0):r_u(x)<r\}.
	\end{align*}			
	The proof is complete.
\end{proof}

\textbf{Acknowledgement.}  The corresponding author C.-L. Xiang is supported by NSFC (No. 12271296) and the NSF of Hubei province, P.R. China (No. 2024AFA061).  G.-F. Zheng is supported by the NSFC (No. 12271195). This work is also  (partly) supported by the Open Research Fund of Hubei Key Laboratory of Mathematical Sciences (No.MPL2025ORG018), Central China Normal University, Wuhan 430079, P. R.China, and supported by the Jiangsu Provincial Scientific Research Center of Applied Mathematics under Grant No. BK20233002.

\end{document}